\documentclass[12pt,a4paper]{article}
\usepackage[left=25mm,right=25mm,top=25mm,bottom=25mm]{geometry}
\usepackage{lmodern}
\usepackage{amsmath,amssymb,amsthm,mathtools}
\usepackage{booktabs,array,longtable}
\usepackage{cite}
\usepackage[colorlinks=true,linkcolor=blue,citecolor=blue,urlcolor=blue]{hyperref}
\usepackage{enumitem}
\usepackage{needspace}
\usepackage{tikz}

\newtheorem{theorem}{Theorem}[section]
\newtheorem{lemma}[theorem]{Lemma}
\newtheorem{proposition}[theorem]{Proposition}

\newtheorem{conjecture}[theorem]{Conjecture}
\theoremstyle{definition}

\newcommand{\N}{N}
\newcommand{\T}{T}
\newcommand{\abs}[1]{\lvert#1\rvert}

\numberwithin{equation}{section}
\numberwithin{figure}{section}
\makeatletter
\newcommand{\preservedpostskip}[1]{\@topsepadd=#1\relax}
\makeatother
\title{Thin Edges in Braces}
\hypersetup{pdftitle={Thin Edges in Braces},pdfauthor={Tong Zhang, Wei Li},pdfsubject={Six-section manuscript with the user-revised opening and an independent cubic-brace theorem}}
\author{Tong Zhang, Wei Li\thanks{ Corresponding author.
		E-mail addresses: liw@nwpu.edu.cn (Wei Li)}\\[0.6em]
\normalsize School of Mathematics and Statistics,\\[0.35em]
\normalsize Northwestern Polytechnical University, Xi'an, Shaanxi, 710129,
P.~R.~China.}
\date{}
\begin{document}
\maketitle
\begin{abstract}
A connected graph with at least two vertices is \emph{matching covered} if each of its
edges is contained in a perfect matching. A bipartite matching covered graph is a
\emph{brace} if no nontrivial cut meets every perfect matching in exactly one edge.
An edge $e$ of a brace $G$ is \emph{thin} if $G-e$ is matching covered and its retract is a brace.
Every brace can be reduced to $K_2$ or $C_4$ by successively deleting thin edges and taking retracts.
Indeed, Carvalho, Lucchesi and Murty proved that every brace of order at least six has at least two thin edges. They
conjectured that there exists a positive constant $c$ such that every brace of order
$n\ge6$ has at least $cn$ thin edges. In this paper, we prove that every simple brace
of order $n\ge6$ has at least $\left\lceil\frac{n+2}{6}\right\rceil$ thin edges and that every brace of order $n$, with
parallel edges allowed, has at least $\frac{n}{13}$ thin edges.
This proves the conjecture.
Moreover, we prove that every cubic brace of order $n$ has at least
$\left\lceil\frac{n+1}{5}\right\rceil$ thin edges for $n\ge6$ and that the coefficient $\frac15$ is best possible.
\end{abstract}



\section{Introduction}


Graphs considered in this paper are loopless, but may have parallel edges.
The \emph{parallel class} of an edge $e$
consists of all edges with the same endvertices as $e$. We denote the vertex set
and edge set of a graph $G$ by $V(G)$ and $E(G)$, respectively.
The \emph{order} of $G$ is $\abs{V(G)}$, and its \emph{size} is $\abs{E(G)}$.
For $v\in V(G)$ and $X\subseteq V(G)$, let $\N_G(v)$ denote the vertices
adjacent to $v$ and $\N_G(X)$ the vertices outside $X$ adjacent to a vertex in $X$.
The \emph{degree} $d_G(v)$ is the number of edges in $G$ incident with the vertex $v$.
For $R\subseteq E(G)$, let $d_R(v)$ be the number of edges of $R$
incident with $v$. A vertex of degree three is \emph{cubic}; a graph
is \emph{cubic} if every vertex is cubic. A vertex set is \emph{stable}
if no two of its vertices are adjacent.
The \emph{induced subgraph} $G[X]$ has
vertex set $X$ and all edges of $G$ between vertices in $X$.
A subgraph is \emph{spanning} if its vertex set is $V(G)$.
The graph $G-X$ is obtained from $G$ by deleting $X$ and all incident edges.
A \emph{path} is a simple graph with at least one vertex whose vertices
can be ordered so that two vertices are adjacent if and only if they are consecutive.
A \emph{cycle} is a connected graph in which every vertex has degree two.
A \emph{forest} is a graph containing no cycle as a subgraph.
A \emph{matching} is a set of edges no two of which have a common endvertex;
it is \emph{perfect} if every vertex is incident with exactly one of its edges.
A graph is \emph{matching covered} if it is connected, has at least two vertices
and every edge belongs to a perfect matching. For undefined notations and terminologies,
please follow \cite{BM}.

For disjoint sets $X,Y\subseteq V(G)$, let $E_G(X,Y)$ consist of
the edges joining $X$ to $Y$, and write
$\partial_G(X)=E_G(X,V(G)\setminus X)$. For a nonempty proper
subset $X$ of $V(G)$, this edge set is a \emph{cut}. If $G$ is
connected, the sets $X$ and $V(G)\setminus X$ are called \emph{shores}.
The cut is \emph{nontrivial} if both $X$ and $V(G)\setminus X$
contain at least two vertices.
A cut of a matching covered graph is \emph{tight} if it meets every
perfect matching in exactly one edge. A matching covered graph with
no nontrivial tight cuts is a \emph{brick} if it is nonbipartite and a
\emph{brace} if it is bipartite. Lov\'asz showed that any two tight cut
decompositions of a matching covered graph yield the same list of
bricks and braces, up to multiple edges \cite{Lovasz}. A brace is
\emph{minimal} if deleting any edge
leaves a graph that is not a brace. An edge $e$ of a matching covered
graph $G$ is \emph{removable} if $G-e$ is matching covered, where
$G-e$ is the graph obtained by deleting $e$ from $G$. The
\emph{bicontraction} of a vertex of
degree two with two distinct
neighbours identifies the three vertices as a single vertex,
removes the resulting loops and retains parallel edges.
The \emph{retract} of a matching covered graph is obtained by
repeatedly bicontracting vertices of degree two until the graph has
order two or has no vertex of degree two. An edge $e$ of a brace
$G$ is \emph{thin} if it is removable and
the retract of $G-e$ is a brace. This retract need not be simple,
even if $G$ is simple \cite{CLM}. We denote \emph{the set of thin edges}
by $\T(G)$ and call an edge \emph{nonthin} if it is not thin.

Lov\'asz's characterization of tight cuts in bipartite matching covered
graphs implies that every tight cut decomposition of such a graph
consists only of braces \cite{Lovasz}. Thus, braces are the bipartite
building blocks in tight cut decompositions of matching covered graphs.
McCuaig proved that every simple brace of order at least six can be
reduced to a biwheel, a bipartite prism or a bipartite M\"obius ladder
by repeatedly deleting a suitable edge and taking the retract \cite{McCuaig}.
Using this theorem, Carvalho, Lucchesi and Murty showed that every
brace on $n$ vertices has at least $\frac{(n-2)^2}{8}$ perfect matchings
\cite{CLMmatchings}.
Lou proved that every minimal simple brace of order at least six has
at least six vertices of degree three \cite{Lou}.
Fabres, Kothari and Carvalho proved that a
minimal brace of order $2n$ has at most $5n-10$ edges for $n\ge6$,
and characterized the graphs for which this bound is attained
\cite{FKC}. Carvalho, Lucchesi and Murty proved
that every brace of order at least six has at least two thin edges.
Although this result guarantees the existence of thin edges, the
lower bound does not increase with the order of the graph.
This gap is particularly important because thin edges are essential building blocks for the inductive generation and structural decomposition of braces \cite{McCuaigPermanent}. A larger number of thin edges facilitates more refined structural analyses and more efficient recursive reductions. Motivated by this significance, they proposed the following conjecture on the number of thin edges in a
brace \cite{CLM}.

\begin{conjecture}\label{conj:linear-thin}\cite{CLM}
There exists a positive constant $c$ such that every brace of order $n\ge6$
has at least $cn$ thin edges.
\end{conjecture}

\begin{figure}[!t]
\centering
\begin{tikzpicture}[line width=0.7pt,
  vertex/.style={circle,draw,inner sep=0pt,minimum size=4.2pt},
  every label/.style={font=\normalsize}]
  \foreach \panel/\offset in {G_1/0,G_2/6} {
    \begin{scope}[xshift=\offset cm]
      \coordinate (a1) at (0,1.6);
      \coordinate (a2) at (0,0);
      \coordinate (a3) at (0,-1.6);
      \coordinate (b1) at (3.2,1.6);
      \coordinate (b2) at (3.2,0);
      \coordinate (b3) at (3.2,-1.6);
      \foreach \j in {2,3} {
        \draw (a1)--(b\j);
      }
      \foreach \i in {2,3} {
        \foreach \j in {1,2,3} {
          \draw (a\i)--(b\j);
        }
      }
      \draw[line width=1.2pt] (a1)--(b1)
        node[midway,above=2pt] {$e$};
      \ifnum\offset=6\relax
        \draw (a1) to[bend left=12] (b2);
        \draw (a2) to[bend left=12] (b1);
      \fi
      \node[vertex,fill=white,label=left:$a_1$] at (a1) {};
      \node[vertex,fill=white,label=left:$a_2$] at (a2) {};
      \node[vertex,fill=white,label=left:$a_3$] at (a3) {};
      \node[vertex,fill=black,label=right:$b_1$] at (b1) {};
      \node[vertex,fill=black,label=right:$b_2$] at (b2) {};
      \node[vertex,fill=black,label=right:$b_3$] at (b3) {};
      \node at (1.6,-2.35) {$\panel$};
    \end{scope}
  }
\end{tikzpicture}
\caption{The edge $e$ is thin in $G_1=K_{3,3}$, but not in $G_2$.}
\label{fig:parallel-thin}
\end{figure}
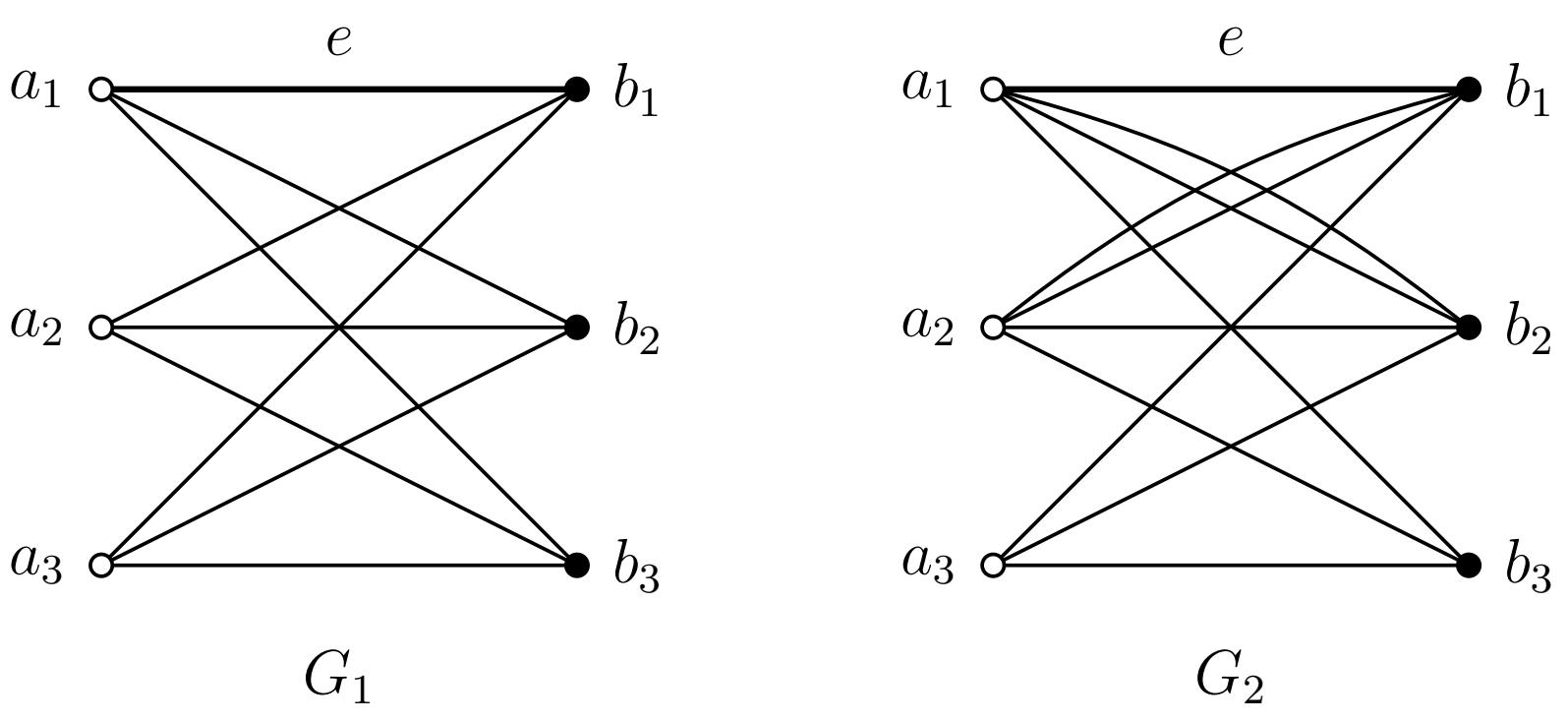

He and Lu proved that every cubic
brace of order $n\ge6$ has a matching of thin edges of size at least
$\frac{n}{10}$, therefore proved the conjecture for cubic braces
\cite{HLcubic}. They later showed that a simple brace of order
$n\ge6$ with $n_3$ vertices of degree three has at least
$n-\frac{5n_3}{2}+1$ thin edges, provided that $0<n_3<\frac{2n}{5}$
\cite{HLforest}. These results do not establish
Conjecture~\ref{conj:linear-thin} for all braces.
In this paper, we prove Conjecture~\ref{conj:linear-thin} by establishing
linear bounds for simple braces and braces with multiple edges, respectively. Our main results are as below.

\begin{theorem}\label{thm:main}
Every simple brace $G$ of order $n\ge6$ satisfies
\[
 \abs{\T(G)}\ge\left\lceil\frac{n+2}{6}\right\rceil.
\]
\end{theorem}

A thin edge of a simple brace need not remain thin after parallel
copies of other edges are added. In Figure~\ref{fig:parallel-thin},
the edge $e=a_1b_1$ is thin in $G_1=K_{3,3}$, but not in the brace $G_2$
obtained by adding a parallel copy of each of $a_1b_2$ and $a_2b_1$.
Thus, Theorem~\ref{thm:main} does not
directly yield a linear lower bound for all braces. In this paper, we obtain
the following lower bound on the number of thin edges in braces with parallel
edges allowed.

\begin{theorem}\label{thm:general}
Every brace $G$ of order $n\ge6$, with parallel edges allowed,
satisfies
\[
 \abs{\T(G)}\ge\left\lceil\frac{n}{13}\right\rceil.
\]
Parallel copies are counted as distinct edges.
\end{theorem}

For cubic braces, we obtain the following stronger lower bound on the number
of thin edges. Moreover, we show that the coefficient $\frac15$ is best possible.

\begin{theorem}\label{cub:thm:main}
Every cubic brace of order $n\ge6$ satisfies
\[
 |T(G)|\ge\left\lceil\frac{n+1}{5}\right\rceil,  \,\,and\,\,
 \qquad\inf_{\substack{G\text{ cubic brace}\\|V(G)|\ge6}}
       \frac{|T(G)|}{|V(G)|}=\frac15.
\]
\end{theorem}

The remainder of this paper is organized as follows.  Section~\ref{sec:prelim} introduces
the preliminary definitions and results. In Section~\ref{ct:section}, we investigate
nonthin edges and critical sets in simple braces, while Section~\ref{sec:short-critical} develops the properties of thin edges and critical sets. Finally, Section~\ref{sec:lower-bounds} provides the proofs of Theorems~\ref{thm:main} and~\ref{thm:general}, and Section~\ref{sec:cubic} gives the proof of Theorem~\ref{cub:thm:main}.


\section{Preliminaries}\label{fd:sec:certificates}\label{sec:prelim}

This section collects the definitions and results used in the proofs.

Let $G$ be a simple brace with bipartition $(A,B)$,
where $\abs{A}=\abs{B}$.
For an odd set $X\subseteq V(G)$, the two sets $X\cap A$ and $X\cap B$
have distinct cardinalities. The larger set is called the \emph{majority part}
of $X$, and the other its \emph{minority part}; we denote these parts by
$X_+$ and $X_-$, respectively. Proposition~\ref{prop:bipartite-tight} gives
a characterization of tight cuts in bipartite matching covered graphs.

\begin{proposition}\label{prop:bipartite-tight}\cite{Lovasz}
Let $G$ be a bipartite matching covered graph and let $C=\partial(X)$
be a cut of $G$, where $\abs{X}$ is odd. Then $C$ is tight if and only if
the majority part of $X$ has one more vertex than its minority part and
every edge of $\partial(X)$ is incident with a vertex in the majority part of $X$.
\end{proposition}

The following lemma describes the perfect matchings and tight cuts of a matching
covered graph and its retract.
\begin{lemma}\label{lem:retract-matchings}
Let $G$ be a matching covered graph, and let $G_1$ be its retract.
Then the following statements hold.
\begin{enumerate}[label=\textup{(\roman*)},leftmargin=*]
\item The restriction of every perfect matching of $G$ to $E(G_1)$ is a perfect
matching of $G_1$, and every perfect matching of $G_1$ can be extended to
a perfect matching of $G$.
\item Let $C_1=\partial_{G_1}(X_1)$ be a cut of $G_1$, and let
$X\subseteq V(G)$ be obtained from $X_1$ by replacing each contraction
vertex in $X_1$ with the set of vertices of $G$ contracted to it.
Then $C_1$ is tight if and only if $\partial_G(X)$ is tight.
\end{enumerate}
\end{lemma}
\begin{proof}
We first prove (i). If $G_1=G$, then the assertion is immediate.
Otherwise, let $G_2$ be obtained from $G$ by bicontracting a vertex $x$
of degree two, with distinct neighbours $y,z$, and let $w$ denote the
contraction vertex. Let $M$ be any perfect matching of $G$. Since $x$ has
degree two, $M$ contains exactly one of $xy$ and $xz$. Adjust notation
so that $xy\in M$. The edge of $M$ incident with $z$ joins $z$ to a
vertex outside $\{x,y,z\}$. Hence $M\setminus\{xy\}$ is a perfect
matching of $G_2$. It follows that the restriction of every perfect
matching of $G$ to $E(G_2)$ is a perfect matching of $G_2$.

Conversely, let $M'$ be any perfect matching of $G_2$, and let $f$ be its
unique edge incident with $w$. In $G$, the edge $f$ joins one of $y,z$ to
a vertex outside $\{x,y,z\}$. Adjust notation so that $f$ is incident with
$y$ in $G$. Then $M'\cup\{xz\}$ is a perfect matching of $G$. Since $G_1$
is obtained from $G$ by repeatedly bicontracting vertices of degree two,
it follows that (i) holds.

We now prove (ii). Set $C=\partial_G(X)$. By the definition of $X$, $C$
and $C_1$ consist of the same edges. Let $M$ be any perfect matching of
$G$, and put $M_1=M\cap E(G_1)$. By (i), $M_1$ is a perfect matching of
$G_1$, and
\[
 |M\cap C|=|M_1\cap C_1|.
\]
By (i), every perfect matching of $G_1$ is the restriction of a perfect
matching of $G$. Hence $C$ is tight if and only if $C_1$ is tight.
\end{proof}

Next, we state a result concerning removable edges in braces.
\begin{theorem}\label{thm:brace-removable}\cite{CLM}
In a brace $G$ on six or more vertices, every edge is removable.
\end{theorem}

The following characterization of braces will be used throughout the proofs.
\begin{theorem}\label{thm:brace-hall}\cite{Lovasz}
A bipartite matching covered graph $G$ with bipartition $(A,B)$, where $\abs{A}=\abs{B}=q$
and $q\ge3$, is a brace if and only if
\[
 \abs{\N(X)}\ge\abs{X}+2,
 \qquad\text{for all }X\subseteq A\text{ such that }0<\abs{X}<q-1.
\]
\end{theorem}

The following two lemmas are consequences of Theorem~\ref{thm:brace-hall}.
\begin{lemma}\label{fd:lem:uncrossing}
Let $G$ be a simple brace with bipartition $(A,B)$ such that
$\abs{A}=\abs{B}=q$ and $q\ge6$. Suppose $S,T$ are nonempty subsets of
the same part of the bipartition, $\abs{N(S)}=\abs{S}+2$ and $\abs{N(T)}=\abs{T}+2$.
If $S\cap T\ne\varnothing$ and $\abs{S\cup T}\le q-2$, then
$\abs{N(S\cap T)}=\abs{S\cap T}+2$, $\abs{N(S\cup T)}=\abs{S\cup T}+2$
and $N(S)\cap N(T)=N(S\cap T)$.
\end{lemma}
\begin{proof}
Since $N(S\cup T)=N(S)\cup N(T)$ and
$N(S\cap T)\subseteq N(S)\cap N(T)$, by the inclusion-exclusion formula
and Theorem~\ref{thm:brace-hall}, we have
\[
\begin{aligned}
 |N(S)|+|N(T)|&=|N(S)\cup N(T)|+|N(S)\cap N(T)|\\
 &\ge|N(S\cup T)|+|N(S\cap T)|\\
 &\ge|S|+|T|+4.
\end{aligned}
\]
Since
\[
 |N(S)|=|S|+2\quad\text{and}\quad |N(T)|=|T|+2,
\]
equality holds throughout. It follows that $|N(S\cap T)|=|S\cap T|+2$
and $|N(S\cup T)|=|S\cup T|+2$. Moreover, equality in the first
inequality gives $|N(S)\cap N(T)|=|N(S\cap T)|$.
Since $N(S\cap T)\subseteq N(S)\cap N(T)$, we obtain
$N(S\cap T)=N(S)\cap N(T)$.
\end{proof}

\begin{lemma}\label{comp:lem:root-deletion}
Let $G$ be a brace with bipartition $(A,B)$, where $|A|=|B|=q$, and let
$X$ be a subset of $A$ or $B$, with $2\le|X|\le q-2$ and
$|N(X)|=|X|+2$. If $N(v)\cap X=\{z\}$, then
\[
 N(X\setminus\{z\})=N(X)\setminus\{v\},\qquad
 |N(X\setminus\{z\})|=|X\setminus\{z\}|+2.
\]
\end{lemma}
\begin{proof}
Since $N(v)\cap X=\{z\}$, we have $v\in N(X)$ and
$N(X\setminus\{z\})\subseteq N(X)\setminus\{v\}$. Hence
\[
 |N(X\setminus\{z\})|\le |N(X)|-1=|X|+1.
\]
Since $1\le |X\setminus\{z\}|\le q-3$, Theorem~\ref{thm:brace-hall} gives
\[
 |N(X\setminus\{z\})|\ge |X\setminus\{z\}|+2=|X|+1.
\]
Thus $|N(X\setminus\{z\})|=|X\setminus\{z\}|+2$.
As $N(X\setminus\{z\})\subseteq N(X)\setminus\{v\}$ and both sets have
$|X|+1$ vertices, we obtain $N(X\setminus\{z\})=N(X)\setminus\{v\}$.
\end{proof}

Here we recall a connectivity property of braces.
\begin{theorem}\label{thm:brace-connectivity}\cite{Plummer}
Every brace of order at least six is $3$-connected.
\end{theorem}

The following lemma describes the retract obtained after deleting an edge from a
simple brace.
\begin{lemma}\label{lem:edge-retract}
Let $G$ be a simple brace of order at least six, and let $e$ be an edge of
$G$ with end vertices $u$ and $v$. The retract of $G-e$ is obtained by
bicontracting the end vertices of $e$ having degree two in $G-e$, if any.
If both end vertices are cubic in $G$, then the sets
$\{u\}\cup N_{G-e}(u)$ and $\{v\}\cup N_{G-e}(v)$ are disjoint.
\end{lemma}
\begin{proof}
Put $H=G-e$, and let $G_1$ be the retract of $H$. The graph $H$ is
matching covered by Theorem~\ref{thm:brace-removable}. By
Theorem~\ref{thm:brace-connectivity}, every vertex of $G$ has degree at
least three. Thus the vertices of degree two in $H$ are precisely the
end vertices of $e$ that are cubic in $G$. If $H$ has no vertex of degree
two, then $G_1=H$ and the assertion is immediate.

Otherwise, adjust notation so that $d_H(u)=2$, and let $u_1,u_2$ be the
two neighbours of $u$ in $H$. Since $G$ is simple, $u_1$ and $u_2$ are
distinct and neither is $v$. Hence both have degree at least three in
$H$. Moreover, they are not adjacent, since $H$ is bipartite.
The vertex obtained by bicontracting $u$ has degree
\[
 d_H(u_1)+d_H(u_2)-2\ge4.
\]
The degrees of vertices outside $\{u,u_1,u_2\}$ do not change, since
parallel edges are retained. If $d_H(v)\ge3$, then every vertex of the
resulting graph has degree at least three, so this graph is $G_1$.

Suppose now that $d_H(v)=2$, and let $v_1,v_2$ be the two neighbours of
$v$ in $H$. Since $u$ and $v$ lie in different parts of the bipartition of $G$
and are not adjacent in $H$, the sets $\{u,u_1,u_2\}$ and $\{v,v_1,v_2\}$ are
disjoint. Thus, after bicontracting $u$, the vertex $v$ still has the two
distinct neighbours $v_1,v_2$, of degrees $d_H(v_1)\ge3$ and
$d_H(v_2)\ge3$, respectively. Bicontracting $v$ now produces a vertex
of degree at least $4$, and the degrees of the other remaining vertices
do not change. Hence every vertex of the resulting graph has degree at
least three, so this graph is $G_1$.
\end{proof}


\section{Nonthin Edges and Critical Sets}\label{ct:section}

In this section, we study the structure of nonthin edges in simple braces. These results will be used to prove Theorem~\ref{thm:main} in Section~\ref{sec:lower-bounds}.

Let $G$ be a simple brace with bipartition $(A,B)$, where $|A|=|B|=q$ and $q\ge6$.
For an edge $uv$ with $u\in A$ and $v\in B$, a set $S\subseteq A$ is a
\emph{nontrivial critical set} with \emph{root} $u$ and \emph{designated neighbour} $v$ if
\begin{equation}\label{eq:route-3-1}
 u\in S,\qquad 2\le|S|\le q-3,\qquad |N(S)|=|S|+2,\qquad N(v)\cap S=\{u\}.
\end{equation}
If $u\in B$ and $v\in A$, conditions~\eqref{eq:route-3-1} define a nontrivial critical set
$S\subseteq B$ with root $u$ and designated neighbour $v$.
The same vertex set may be a nontrivial critical set with different choices of root and designated neighbour.

\subsection{A Characterization of Nonthin Edges}
We begin with a characterization of nonthin edges in terms of nontrivial critical sets.

\begin{proposition}\label{res:route-3-1}
Let $G$ be a simple brace with bipartition $(A,B)$, where $|A|=|B|=q$ and $q\ge6$.
Let $uv\in E(G)$, with $u\in A$ and $v\in B$. The edge $uv$ is nonthin if
and only if there are $S\subseteq A$ and $R\subseteq B$ such that
$|S|,|R|\ge2$, $|S|+|R|=q-1$ and $E_G(S,R)=\{uv\}$.
\end{proposition}
\begin{proof}
Put $H=G-uv$ and let $G_1$ be its retract. Both are bipartite matching
covered graphs by Theorem~\ref{thm:brace-removable} and Lemma~\ref{lem:retract-matchings}(i). By Lemma~\ref{lem:edge-retract}, the
bicontractions identify $K_x=\{x\}\cup N_H(x)$ for each cubic end vertex $x$
of $uv$; if both end vertices are cubic, $K_u,K_v$ are disjoint.

Suppose $uv$ is nonthin. Take a nontrivial tight cut of $G_1$ and replace
each contraction vertex by its set $K_x$. The resulting shores $X,Y$
define a tight cut of $H$ by Lemma~\ref{lem:retract-matchings}(ii). Their majority parts exceed their
minority parts by one, and all cut edges join the majority parts, by
Proposition~\ref{prop:bipartite-tight}. Those majority parts lie in different parts of the
bipartition since $|A|=|B|$. If $u,v$ were in one shore or in the two
majority parts, adding $uv$ would leave the cut tight in $G$, contrary
to $G$ being a brace. An edge between different shores cannot join a
majority part to a minority part. Thus, after interchanging $X,Y$,
$u\in X\cap A$ and $v\in Y\cap B$ are in the minority parts.

Put $S=X\cap A$ and $R=Y\cap B$. Then
$E_H(S,R)=\varnothing$ and $|B\setminus R|=|S|+1$, so
$E_G(S,R)=\{uv\}$ and $|S|+|R|=q-1$.
If $|S|=1$, then $S=\{u\}$ and $|X\cap B|=2$.
All neighbours of $u$ in $H$ lie in $X\cap B$; simplicity and minimum
degree three in $G$ give $d_H(u)=2$ and $X=K_u$.
Its image would be a singleton shore in $G_1$, a contradiction.
Hence $|S|\ge2$, and the same argument gives $|R|\ge2$.

Conversely, suppose $S,R$ satisfy the stated conditions. In $H$, put
$X=S\cup(B\setminus R)$ and $Y=R\cup(A\setminus S)$.
The equality $|S|+|R|=q-1$ makes $B\setminus R$ and $A\setminus S$
the majority parts of their shores, each exceeding its minority part
by one. Since $E_H(S,R)=\varnothing$, Proposition~\ref{prop:bipartite-tight} makes this cut
tight. Each defined $K_u$ lies in $X$ and each defined $K_v$ in $Y$.
Thus the bicontractions give a tight cut of $G_1$ by Lemma~\ref{lem:retract-matchings}(ii).
Each shore has at least five vertices before the bicontractions and loses at most two vertices,
so the resulting cut is nontrivial. Hence $uv$ is nonthin.
\end{proof}

We say that $S\subseteq A$ and $R\subseteq B$ are \emph{paired} for $uv$ if
$|S|,|R|\ge2$, $|S|+|R|=q-1$, and $E_G(S,R)=\{uv\}$.
The following lemma gives some properties of paired sets.
\begin{lemma}\label{res:route-3-2}
Let $G$ be a simple brace with bipartition $(A,B)$, where $|A|=|B|=q$ and $q\ge6$.
Let $uv\in E(G)$, with $u\in A$ and $v\in B$. Then the following statements hold.
\begin{enumerate}[label=\textup{(\roman*)},leftmargin=*]
\item Paired sets $S,R$ for $uv$ are nontrivial critical sets with respective roots
$u,v$ and respective designated neighbours $v,u$, and
$N(S)=(B\setminus R)\cup\{v\}$ and $N(R)=(A\setminus S)\cup\{u\}$.
\item Every nontrivial critical set with root $u$ and designated neighbour $v$
has a unique paired nontrivial critical set with root $v$ and designated neighbour $u$, and conversely.
\item If $S_1,R_1$ and $S_2,R_2$ are paired for the same edge $uv$, then
$S_1\subseteq S_2$ if and only if $R_2\subseteq R_1$.
\end{enumerate}
\end{lemma}
\begin{proof}
We first prove (i). Since $|S|,|R|\ge2$ and $|S|+|R|=q-1$, we have
$2\le|S|,|R|\le q-3$. Since $E_G(S,R)=\{uv\}$, we have $u\in S$,
$v\in R$, $N(v)\cap S=\{u\}$ and $N(u)\cap R=\{v\}$.
Moreover, $N(S)\subseteq(B\setminus R)\cup\{v\}$, so
$|N(S)|\le q-|R|+1=|S|+2$. By Theorem~\ref{thm:brace-hall},
$|N(S)|\ge|S|+2$. Hence $|N(S)|=|S|+2$ and
$N(S)=(B\setminus R)\cup\{v\}$.
Similarly, $|N(R)|=|R|+2$ and $N(R)=(A\setminus S)\cup\{u\}$.
Thus $S$ and $R$ are nontrivial critical sets with respective roots $u,v$
and respective designated neighbours $v,u$.

For (ii), let $S$ be a nontrivial critical set with root $u$ and designated neighbour $v$, and put
$R=(B\setminus N(S))\cup\{v\}$. Then $|R|=q-|S|-1\ge2$ and
$|S|+|R|=q-1$. Since $N(v)\cap S=\{u\}$, we also have
$E_G(S,R)=\{uv\}$. Thus $S$ and $R$ are paired, and (i) implies that
$R$ is a nontrivial critical set with root $v$ and designated neighbour $u$.
Moreover, by (i), any set paired with $S$ for $uv$ must equal
$(B\setminus N(S))\cup\{v\}$. Hence $R$ is unique.
By symmetry, every nontrivial critical set with root $v$ and designated neighbour $u$
is paired with a unique nontrivial critical set with root $u$ and designated neighbour $v$.

To prove (iii), note that (i) gives $R_i=(B\setminus N(S_i))\cup\{v\}$ and
$S_i=(A\setminus N(R_i))\cup\{u\}$ for $i=1,2$.
If $S_1\subseteq S_2$, then $N(S_1)\subseteq N(S_2)$, and hence $R_2\subseteq R_1$.
Conversely, if $R_2\subseteq R_1$, then $N(R_2)\subseteq N(R_1)$,
and hence $S_1\subseteq S_2$.
\end{proof}

The following proposition is a consequence of Proposition~\ref{res:route-3-1}.

\begin{proposition}\label{res:route-3-4}
If $H\subseteq G$ are simple braces on the same set of at least six vertices,
then $T(H)\subseteq T(G)$.
\end{proposition}
\begin{proof}
Let $(A,B)$ be a bipartition of $G$, where $|A|=|B|=q$ and $q\ge3$.
Since $H$ is a spanning subgraph of $G$, $(A,B)$ is also a bipartition of $H$.
The proof of Proposition~\ref{res:route-3-1} applies also when $q\ge3$.
Let $e\in T(H)$, and suppose, to the contrary, that $e$ is nonthin in $G$.
By Proposition~\ref{res:route-3-1}, there are $S\subseteq A$ and $R\subseteq B$
such that $|S|,|R|\ge2$, $|S|+|R|=q-1$ and $E_G(S,R)=\{e\}$.
Since $H\subseteq G$ and $e\in E(H)$, we have $E_H(S,R)=\{e\}$.
Proposition~\ref{res:route-3-1}, applied to $H$, implies that $e$ is nonthin in $H$,
contrary to $e\in T(H)$. Hence $e\in T(G)$, and therefore $T(H)\subseteq T(G)$.
\end{proof}

The following lemma allows us to
restrict the later arguments to minimal simple braces.

\begin{lemma}\label{prop:spanning-minimal-brace}
Every simple brace $G$ has a spanning subgraph $H$ that is a minimal brace and satisfies
$|T(G)|\ge|T(H)|$.
\end{lemma}
\begin{proof}
If $G$ is isomorphic to $K_2$ or $C_4$, take $H=G$.
Otherwise, choose a spanning brace $H\subseteq G$ with as few
edges as possible. It is minimal, and Proposition~\ref{res:route-3-4} gives
$|T(G)|\ge|T(H)|$.
\end{proof}

\subsection{Properties of Critical Sets}
We begin with some notation. Throughout this subsection, let $G$ be a minimal simple brace with
bipartition $(A,B)$, where $|A|=|B|=q$ and $q\ge6$.

Let $C$ be the set of cubic vertices incident with thin edges and
$V_{\ge4}=\{x\in V(G):d(x)\ge4\}$. For each cubic vertex $x\notin C$, write
$N(x)=\{y_1,y_2,y_3\}$. All three edges incident with $x$ are nonthin, so
Proposition~\ref{res:route-3-1} and Lemma~\ref{res:route-3-2} allow us to choose,
for $i=1,2,3$, a nontrivial critical set $S_i$ with root $y_i$ and designated neighbour $x$.
Let $D$ consist of those $x$ for which $S_1,S_2,S_3$ are pairwise disjoint, and $P$ of the remaining
cubic vertices outside $C$. Put $U=V_{\ge4}\cup D$ and $p=|P|$.
The next lemma shows that $D,P$ do not depend on the choices.

\begin{lemma}\label{res:route-3-6}
The following statements hold.
\begin{enumerate}[label=\textup{(\roman*)},leftmargin=*]
\item Let $x$ be a vertex of $G$ with $d(x)\ge4$, and let $xy_1,xy_2$ be distinct nonthin edges.
For $i=1,2$, let $S_i$ be a nontrivial critical set with root $y_i$ and designated neighbour $x$.
Then $S_1\cap S_2=\varnothing$.
\item The sets $D$ and $P$ do not depend on the choices in their definition.
\end{enumerate}
\end{lemma}
\begin{proof}
Let $x\in V(G)$, and let $W$ be the part of the bipartition containing $N(x)$.
Let $xy_1,xy_2$ be distinct nonthin edges and, for $i=1,2$, let $S_i$ be a
nontrivial critical set with root $y_i$ and designated neighbour $x$. Suppose that $S_1\cap S_2\ne\varnothing$.
Since $N(x)\cap S_i=\{y_i\}$ and $y_1\ne y_2$, the vertex $x$ has no neighbour
in $S_1\cap S_2$. If $|S_1\cup S_2|\le q-2$, Lemma~\ref{fd:lem:uncrossing} gives
$x\in N(S_1)\cap N(S_2)=N(S_1\cap S_2)$, a contradiction.
Also, $S_1\cup S_2$ contains no vertex of $N(x)\setminus\{y_1,y_2\}$, so
$q-1\le|S_1\cup S_2|\le q-d(x)+2$. Since $d(x)\ge3$, we obtain $d(x)=3$ and
$S_1\cup S_2=W\setminus\{y_3\}$, where $N(x)=\{y_1,y_2,y_3\}$.
Thus, if $d(x)\ge4$, the sets $S_1,S_2$ cannot intersect. This proves (i).

For (ii), suppose that $x$ is cubic with no incident thin edge. Write
$N(x)=\{y_1,y_2,y_3\}$ and, for $i=1,2,3$, choose a nontrivial critical set
$S_i$ with root $y_i$ and designated neighbour $x$. Suppose that $S_1\cap S_2\ne\varnothing$.
Since $|S_3|\ge2$ and $y_3\in S_3$, the set $S_3$ intersects
$S_1\cup S_2=W\setminus\{y_3\}$. Interchanging the subscripts $1,2$ if necessary,
assume $S_1\cap S_3\ne\varnothing$. The preceding argument gives
$S_1\cup S_3=W\setminus\{y_2\}$. Since $|S_1\cup N(x)|=|S_1|+2\le q-1$,
there is a vertex $z\in W\setminus(S_1\cup N(x))$.
The two union equalities and $z\notin S_1$ imply $z\in S_2\cap S_3$.
Thus the three sets intersect pairwise, and $S_i\cup S_j=W\setminus\{y_k\}$
for $\{i,j,k\}=\{1,2,3\}$.

Let $S'_1$ be any nontrivial critical set
with root $y_1$ and designated neighbour $x$. If $S_2\cap S_3\ne\varnothing$, the preceding argument shows that
$S'_1,S_2,S_3$ intersect pairwise. If $S_2\cap S_3=\varnothing$, then $S'_1$
cannot intersect either set, since such an intersection would imply
$S_2\cap S_3\ne\varnothing$. Hence replacing $S_1$ does not change whether $x$
belongs to $D$ or $P$. Replacing $S_2$ and $S_3$ in turn proves that $D,P$ do not depend on these choices.
This proves (ii).
\end{proof}

\begin{lemma}\label{res:route-3-7}
Let $x\in P$ with $N(x)=\{y_1,y_2,y_3\}$. For $i=1,2,3$, let $S_i$ be a nontrivial
critical set with root $y_i$ and designated neighbour $x$.
If $W$ is the part of the bipartition containing $N(x)$,
then $S_i\cup S_j=W\setminus\{y_k\}$ for
$\{i,j,k\}=\{1,2,3\}$.
\end{lemma}
\begin{proof}
Let $\{i,j,k\}=\{1,2,3\}$. Since $x\in P$, the proof of
Lemma~\ref{res:route-3-6}(ii) shows that $S_i\cap S_j\ne\varnothing$.
Since $N(x)\cap S_i=\{y_i\}$ and
$N(x)\cap S_j=\{y_j\}$, the vertex $x$ has no neighbour in $S_i\cap S_j$.
If $|S_i\cup S_j|\le q-2$, Lemma~\ref{fd:lem:uncrossing} gives
$x\in N(S_i)\cap N(S_j)=N(S_i\cap S_j)$, a contradiction.
Thus $|S_i\cup S_j|\ge q-1$. The same neighbourhood conditions give $y_k\notin S_i\cup S_j$.
Since $S_i\cup S_j\subseteq W$ and $|W|=q$, it follows that $S_i\cup S_j=W\setminus\{y_k\}$.
\end{proof}

The next lemma shows that $G[U]$ is a forest.

\begin{lemma}\label{res:route-3-8}
The following statements hold.
\begin{enumerate}[label=\textup{(\roman*)},leftmargin=*]
\item Every thin edge of $G$ has a cubic end vertex.
\item Let $x\in U$, and let $xy_1,xy_2$ be distinct nonthin edges.
For $i=1,2$, let $S_i$ be a nontrivial critical set with root $y_i$ and designated neighbour $x$.
Then $S_1\cap S_2=\varnothing$ and $|S_1|+|S_2|\le q-2$.
\item The graph $G[U]$ is a forest.
\end{enumerate}
\end{lemma}
\begin{proof}
We first prove (i). Suppose that a thin edge $xy$ has no cubic end vertex.
Then every vertex of $G-xy$ has degree at least three, so $G-xy$ is its own retract.
Since $xy$ is thin, $G-xy$ is a brace, contrary to the minimality of $G$.

To prove (ii), let $W$ be the part of the bipartition containing $N(x)$.
If $x\in V_{\ge4}$, Lemma~\ref{res:route-3-6}(i) gives $S_1\cap S_2=\varnothing$.
Since $N(x)\cap S_i=\{y_i\}$ for $i=1,2$, at least two neighbours of $x$ lie outside $S_1\cup S_2$.
If $x\in D$, write $N(x)=\{y_1,y_2,y_3\}$ and choose a nontrivial critical set
$S_3$ with root $y_3$ and designated neighbour $x$.
By Lemma~\ref{res:route-3-6}(ii) and the definition of $D$, the sets $S_1,S_2,S_3$ are pairwise disjoint.
Since $|S_3|\ge2$, at least two vertices of $W$ lie outside $S_1\cup S_2$.
Thus, in both cases, $S_1\cap S_2=\varnothing$ and $|S_1|+|S_2|\le |W|-2=q-2$.

For (iii), every thin edge has an end vertex in $C$ by (i).
Since $U\cap C=\varnothing$, every edge of $G[U]$ is nonthin.
Suppose that $G[U]$ contains a cycle $x_0x_1\cdots x_{s-1}x_0$.
For $i=0,\ldots,s-1$, choose paired nontrivial critical sets $S_i,R_i$ with respective
roots $x_{i-1},x_i$ and respective designated neighbours $x_i,x_{i-1}$,
where subscripts are taken modulo $s$. By the definition of paired sets, $|S_i|+|R_i|=q-1$.
The sets $S_i,R_{i+1}$ both have designated neighbour $x_i$, so (ii) gives
$|S_i|+|R_{i+1}|\le q-2$. Summing over the vertices of the cycle, we obtain
\[
s(q-1)=\sum_{i=0}^{s-1}(|S_i|+|R_i|)
=\sum_{i=0}^{s-1}(|S_i|+|R_{i+1}|)\le s(q-2),
\]
a contradiction. Hence $G[U]$ is a forest.
\end{proof}

The following lemma is a consequence of Lemma~\ref{res:route-3-8}(ii).

\begin{lemma}\label{res:route-3-9}
Let $x_0\cdots x_s$ be a path whose edges are nonthin, with $s\ge2$ and
$x_1,\ldots,x_{s-1}\in U$. If $F_0,F_s$ are nontrivial critical sets with
respective roots $x_0,x_s$ and respective designated neighbours $x_1,x_{s-1}$,
then $|F_0|+|F_s|+s\le q$.
\end{lemma}
\begin{proof}
For each edge $x_{i-1}x_i$ of the path, choose paired nontrivial critical sets $S_i,R_i$
rooted at $x_{i-1},x_i$, respectively. Take $S_1=F_0$ and $R_s=F_s$.
Since $|S_i|+|R_i|=q-1$ by Proposition~\ref{res:route-3-1}, and $S_i,R_{i+1}$
have roots $x_{i-1},x_{i+1}$, respectively, and both have designated neighbour $x_i$,
Lemma~\ref{res:route-3-8}(ii) gives
\[
 \sum_{i=1}^{s-1}(|S_i|+|R_{i+1}|)
 =|F_0|+|F_s|+(s-2)(q-1)\le(s-1)(q-2).
\]
Rearranging gives
\begin{equation}\label{eq:route-3-2}
 |F_0|+|F_s|+s\le q.
\end{equation}
\end{proof}

We next consider nontrivial critical sets with a common root in $P$.

\begin{lemma}\label{res:route-3-10}\label{res:route-3-11}
Let $u\in P$ with $N(u)=\{y_1,y_2,y_3\}$. For $i=1,2,3$, let $S_i$ be a nontrivial
critical set with root $u$ and designated neighbour $y_i$, and put $Y_i=N(S_i)\setminus N(u)$.
Then the following statements hold.
\begin{enumerate}[label=\textup{(\roman*)},leftmargin=*]
\item $S_i\cap S_j=\{u\}$ for distinct $i,j\in\{1,2,3\}$.
\item The sets $Y_1,Y_2,Y_3$ are pairwise disjoint, and $|Y_i|=|S_i|-1$ for $i=1,2,3$.
\item $\sum_{i=1}^3|S_i|\le q$, and $|S_i|\le q-4$ for $i=1,2,3$.
\end{enumerate}
\end{lemma}
\begin{proof}
Assume $u\in A$ and, for $i=1,2,3$, let $R_i$ be paired with $S_i$.
Lemma~\ref{res:route-3-7} gives
$R_i\cup R_j=B\setminus\{y_k\}$ whenever $\{i,j,k\}=\{1,2,3\}$.

To prove (i), suppose that $z\in(S_i\cap S_j)\setminus\{u\}$ for distinct $i,j\in\{1,2,3\}$.
Since $E_G(S_i,R_i)=\{uy_i\}$ and $E_G(S_j,R_j)=\{uy_j\}$, the vertex $z$ has no neighbour
in $R_i\cup R_j$. Thus $N(z)\subseteq B\setminus(R_i\cup R_j)=\{y_k\}$, contrary to $d(z)\ge3$.
As $u\in S_i\cap S_j$, this proves (i).

For (ii), Lemma~\ref{res:route-3-2}(i) gives $N(S_i)=(B\setminus R_i)\cup\{y_i\}$,
and hence $Y_i=B\setminus(R_i\cup N(u))$.
Since $R_i\cup R_j\cup N(u)=B$ for $i\ne j$, the sets $Y_1,Y_2,Y_3$ are pairwise disjoint.
Also, $N(u)\subseteq N(S_i)$, so $|Y_i|=|N(S_i)|-|N(u)|=|S_i|-1$.

Finally, the sets $Y_1,Y_2,Y_3$ are disjoint subsets of $B\setminus N(u)$.
By (ii), $\sum_{i=1}^3(|S_i|-1)=\sum_{i=1}^3|Y_i|\le q-3$, and hence $\sum_{i=1}^3|S_i|\le q$.
For each $i$, the other two sets have at least two vertices each, so $|S_i|\le q-4$.
This proves (iii).
\end{proof}

We next consider maximal nontrivial critical sets. Let $G$ be a minimal simple brace
with bipartition $(A,B)$, where $|A|=|B|=q$ and $q\ge6$.
We use the notation $P$ introduced at the beginning of this subsection.
Every vertex of $P$ is cubic and is incident with no thin edge.
A nontrivial critical set $S$ is \emph{maximal} if no nontrivial critical set
with the same root and designated neighbour properly contains $S$.
Since $G$ is finite, every nontrivial critical set is contained in a maximal nontrivial critical set
with the same root and designated neighbour.
Put $b_q=\lfloor\frac{q-2}{2}\rfloor$. A nontrivial critical set $S$ is \emph{small} if $|S|\le b_q$,
and \emph{large} otherwise.
\begin{lemma}\label{res:route-3-12}
Let $S$ be a maximal nontrivial critical set with root $u$ and designated neighbour $v$, and let $X$ be a set in the
same part with $|N(X)|=|X|+2$. If $S\cap X\ne\varnothing$, $u\notin X$ and $|S\cup X|\le q-3$, then
$X\subseteq S\setminus\{u\}$.
\end{lemma}
\begin{proof}
By Lemma~\ref{fd:lem:uncrossing}, $|N(S\cup X)|=|S\cup X|+2$ and $N(S)\cap N(X)=N(S\cap X)$.
Since $u\notin X$ and $N(v)\cap S=\{u\}$, the vertex $v$ has no neighbour in $S\cap X$.
As $v\in N(S)$ and $N(S)\cap N(X)=N(S\cap X)$, it follows that $v\notin N(X)$.
Thus $N(v)\cap(S\cup X)=\{u\}$. Since $2\le|S\cup X|\le q-3$,
the set $S\cup X$ is a nontrivial critical set with root $u$ and designated neighbour $v$.
By the maximality of $S$, it follows that $S\cup X=S$.
Hence $X\subseteq S\setminus\{u\}$, since $u\notin X$.
In particular, the conclusion holds for small $S$ and
$|X|\le b_q$, since $|S\cup X|\le2b_q-1\le q-3$.
\end{proof}

\begin{lemma}\label{res:route-3-13}
Let $S$ be a small maximal nontrivial critical set with root $u$ and designated neighbour $v$. A cubic vertex outside
$S$ with at least two neighbours in $N(S)$ has exactly two such neighbours, one of
which is $v$.
\end{lemma}
\begin{proof}
Let $x\in V(G)\setminus S$ be a cubic vertex with at least two neighbours in $N(S)$.
Since $G$ is bipartite, $x$ lies in the same part as the vertices of $S$.
Put $S'=S\cup\{x\}$. We have $|S'|\le b_q+1\le q-3$.
If all three neighbours of $x$ belonged to $N(S)$, then $|N(S')|=|S'|+1$, contrary to
Theorem~\ref{thm:brace-hall}. Thus $x$ has exactly two neighbours in $N(S)$.
Its third neighbour lies outside $N(S)$, so $|N(S')|=|N(S)|+1=|S'|+2$.
Suppose that $xv\notin E(G)$. Then $N(v)\cap S'=\{u\}$, so $S'$ is a nontrivial critical set
with root $u$ and designated neighbour $v$.
Since $S\subsetneq S'$, this contradicts the maximality of $S$. Hence $xv\in E(G)$, as required.
\end{proof}

For each $u\in P$ and each $v\in N(u)$, select one maximal nontrivial critical set
with root $u$ and designated neighbour $v$.

\begin{lemma}\label{res:route-3-17}
Every cubic vertex of $G$ has at most one neighbour in $P$.
\end{lemma}
\begin{proof}
Suppose that a cubic vertex $v$ has distinct neighbours $u,w\in P$,
and write $N(v)=\{u,w,z\}$. Let $S,T$ be the two chosen maximal nontrivial critical sets
with root $u$ and designated neighbours different from $v$.
For each $F\in\{S,T\}$, the vertex $v$ belongs to $N(F)$ and is not the designated
neighbour of $F$. By Lemma~\ref{comp:lem:root-deletion}, $v\in N(F\setminus\{u\})$.
Thus each of $S,T$ contains at least one of $w,z$.
By Lemma~\ref{res:route-3-10}(i), $S\cap T=\{u\}$.
Adjust notation so that $w\in S$ and $z\in T$. In particular, $w\notin T$.
Similarly, let $K',J$ be the two chosen maximal nontrivial critical sets with root $w$
and designated neighbours different from $v$, with $u\in K'$ and $z\in J$.
Since $K'\cap J=\{w\}$ by Lemma~\ref{res:route-3-10}(i), we have $z\notin K'$.

Interchanging $u$ with $w$, $S$ with $K'$, and $T$ with $J$ if necessary, assume $|K'|\le|S|$.
The chosen critical set with root $u$ and designated neighbour $v$ has at least two vertices, so Lemma~\ref{res:route-3-11}(iii) gives
$|S|+|T|\le q-2$. Since $u\in T\cap K'$, we have
\[
 |T\cup K'|\le |T|+|K'|-1\le |T|+|S|-1\le q-3.
\]
As $w\notin T$, Lemma~\ref{res:route-3-12}, applied to $K'$ and $T$, gives
$T\subseteq K'\setminus\{w\}$, contrary to $z\in T$ and $z\notin K'$.
\end{proof}

\subsection{The Number of Small Critical Sets}
The results of this subsection will be used to prove the lower bound on the number of thin edges.
We begin with some definitions.
Let $G$ be a minimal simple brace with bipartition $(A,B)$, where $|A|=|B|=q$ and $q\ge6$.
We use the set $P$ defined in Subsection~3.2.
Every vertex of $P$ is cubic and is incident with no thin edge.

For each $u\in P$ and each $v\in N(u)$, by Proposition~\ref{res:route-3-1},
Lemma~\ref{res:route-3-2}(i) and the definition of maximality, there exists at least one
maximal nontrivial critical set with root $u$ and designated neighbour $v$, since $G$ is finite.
Let $S_{u,v}$ denote the maximal nontrivial critical set
selected with root $u$ and designated neighbour $v$. Let
\[
 \mathcal S=\{(S_{u,v},u,v):u\in P,\ v\in N(u)\}.
\]
Put $p=|P|$ and $\sigma=\sum_{v\in V(G)}(d(v)-3)$.
Since every vertex of $P$ has three neighbours, $|\mathcal S|=3p$.

\begin{lemma}\label{res:route-3-18}
Let $G$ be a minimal simple brace of order $2q$, where $q\ge6$, and suppose $p>0$.
Let $b_0$ denote the number of triples $(S,u,v)\in\mathcal S$ such that $S$ is small and $|S\cap P|\ge2$.
Then $b_0\le p-1+\frac{\sigma}{2}$.
\end{lemma}
\begin{proof}
Let $S,T$ be small maximal nontrivial critical sets with respective roots $u,v$,
where $u\notin T$ and $v\notin S$. Suppose that $S$ and $T$ have a vertex in common.
Then $|S\cup T|\le2b_q-1\le q-3$. Lemma~\ref{res:route-3-12} gives $T\subseteq S\setminus\{u\}$.
Since $v\in T$, this implies $v\in S$, contrary to $v\notin S$. Hence $S\cap T=\varnothing$.

We shall construct vertex sets $K_1,\ldots,K_m$, where $m$ denotes the total number of sets obtained,
and prove \mbox{$b_0\le p-1+m$}. Put $Q=P$.
Assume, inductively, that each set $K_i$ already constructed is contained in a small
maximal nontrivial critical set $S_i$, with root $x_i$ and designated neighbour $z_i$,
such that $(S_i,x_i,z_i)\in\mathcal S$ and $S_i\cap Q=\{x_i\}$.

If there is a triple $(F,u,v_0)\in\mathcal S$ with $u\in Q$, $|F|\le b_q$
and $|F\cap Q|\ge2$, select one for which $F$ has minimum cardinality, and put $X=F\cap Q$.
Let $v\in X\setminus\{u\}$. If a triple $(E,v,b)\in\mathcal S$ satisfies
$|E|\le b_q$ and $|E\cap Q|\ge2$, then $u\in E$.
For otherwise, since $v\in E\cap F$, Lemma~\ref{res:route-3-12} gives
$E\subseteq F\setminus\{u\}$, so $|E|<|F|$, contrary to the choice of $F$.
There is at most one such triple with root $v$.
Indeed, by the definition of $\mathcal S$, two distinct triples with root $v$ have different designated neighbours.
By Lemma~\ref{res:route-3-10}(i), their critical sets have only $v$ in common,
whereas both contain $u\ne v$. If such a triple exists, let $E_v$ denote its critical set.

Suppose that $E_v$ exists for every $v\in X\setminus\{u\}$, and select a set $E_w$
of minimum cardinality.
If $w\notin E_v$, then $u\in E_v\cap E_w$, so Lemma~\ref{res:route-3-12} gives
$E_v\subseteq E_w\setminus\{w\}$, contrary to the choice of $w$.
Thus $K=F\cap\bigcap_{v\in X\setminus\{u\}}E_v$ contains $u,w$.
Since $|F|,|E_v|\le b_q$ and $2b_q\le q-2$, repeated applications of
Lemma~\ref{fd:lem:uncrossing} give $|N(K)|=|K|+2$.
For the next index $j$, let $K_j$ be obtained from $K$ by deleting the vertices of
$(K\cap Q)\setminus\{u,w\}$, one at a time.
For each such vertex $v$, the designated neighbour of $E_v$ has only $v$ as a neighbour in $K$, since $K\subseteq E_v$.
The vertices $u,w$ are not deleted, so Lemma~\ref{comp:lem:root-deletion} gives $|N(K_j)|=|K_j|+2$.
Put $R_j=\{u,w\}$, and let $y_j(u),y_j(w)$ be the designated neighbours of $F,E_w$, respectively.
Then $|K_j|\le b_q$ and $K_j\cap Q=R_j$.
If $E_v$ does not exist for some $v\in X\setminus\{u\}$, no new set $K_j$ is constructed.

Suppose that $K_j$ has been constructed, and let $i<j$.
If $x_i\notin X$, then $x_i\notin F$, since $x_i\in Q$ and $X=F\cap Q$.
Also, $u\notin S_i$, since $u\in Q\setminus\{x_i\}$ and $S_i\cap Q=\{x_i\}$.
The first paragraph gives $S_i\cap F=\varnothing$, so $K_i\cap K_j=\varnothing$.
If $x_i\in X$, the set $S_i$ has root $x_i$ and $S_i\cap Q=\{x_i\}$.
When $x_i=u$, the set $F$ has the same root and $|F\cap Q|\ge2$;
when $x_i\ne u$, the set $E_{x_i}$ has the same root and $|E_{x_i}\cap Q|\ge2$.
In either case, the two critical sets are distinct.
By the definition of $\mathcal S$, they have different designated neighbours.
Lemma~\ref{res:route-3-10}(i) therefore gives $K_i\cap K_j\subseteq\{x_i\}$.

When $K_i$ was constructed, $K_i\cap Q=R_i$, and one vertex of $R_i$ was then deleted from $Q$.
Since vertices are only deleted from $Q$, we still have $K_i\cap Q\subseteq R_i$.
Suppose that $K_i\cap K_j\ne\varnothing$, so $K_i\cap K_j=\{x_i\}$.
As $x_i\in Q$ and $K_j\cap Q=R_j$, it follows that $x_i\in R_i\cap R_j$.
Thus $x_i$ is the vertex of $R_i$ retained in $Q$ when $K_i$ was constructed.
Under the replacement rule given below, the triple associated with $K_i$ has not changed,
since $x_i$ has not been deleted. Hence $z_i=y_i(x_i)$.
If $x_i=u$, the designated neighbour of $F$ is $y_j(x_i)$;
if $x_i=w$, that of $E_w$ is $y_j(x_i)$.
The different designated neighbours obtained above are therefore $y_i(x_i)$ and $y_j(x_i)$.
As $R_i\subseteq K_i$ and $R_j\subseteq K_j$, we conclude that
$K_i\cap K_j=R_i\cap R_j$ and $y_i(x)\ne y_j(x)$ for every $x\in R_i\cap R_j$.

Now delete $X\setminus\{u\}$ from $Q$.
If $x_i$ is deleted, then $x_i\in S_i\cap F$ and $u\notin S_i$.
Lemma~\ref{res:route-3-12} gives $S_i\subseteq F\setminus\{u\}$.
For each such $K_i$, replace $(S_i,x_i,z_i)$ by $(F,u,v_0)$.
If a new set $K_j$ was constructed, put $(S_j,x_j,z_j)=(F,u,v_0)$.
Leave all other triples unchanged.
After the deletion, $F\cap Q=\{u\}$, so $K_i\subseteq S_i$ and $S_i\cap Q=\{x_i\}$ still hold for every $K_i$.
Repeat this procedure until no triple satisfies the conditions for selecting $F$.
Since $Q$ is finite and each step deletes at least one vertex while retaining $u$,
the process ends with $Q\ne\varnothing$.

Consider the bipartite graph with one vertex for each $K_i$ and one vertex for each
$r\in\bigcup_iR_i$, where $K_i$ and $r$ are adjacent if and only if $r\in R_i$.
We show that this graph is a forest.
Suppose that it contains a cycle, and select a root $r$ on the cycle that was deleted from $Q$ first.
Such a root exists because one vertex of each $R_i$ was deleted in the step in which $K_i$ was constructed.
Let $K_i,K_j$ be the two neighbours of $r$ on the cycle.
If $K_i$ was constructed before the step in which $r$ was deleted, the other vertex of $R_i$,
which is also on the cycle, was deleted before $r$, a contradiction.
Nor could $K_i$ have been constructed after that step, since $r$ was no longer in $Q$.
Thus $K_i$ was constructed in that step, and the same holds for $K_j$.
This is impossible, since at most one set $K_i$ is constructed at each step.

At the beginning of each step, count the triples $(S,a,b)\in\mathcal S$ with
$a\in Q$, $|S|\le b_q$ and $|S\cap Q|\ge2$.
Put $k=|X|-1$, and let $t$ be the number of vertices $v\in X\setminus\{u\}$ for which $E_v$ exists.
By the uniqueness proved above, exactly $t$ of the counted triples have their roots among these $k$ vertices.
Among triples whose roots remain in $Q$, only $(F,u,v_0)$ ceases to be counted.
Indeed, let $(J,a,b)$ be a triple counted before but not after the deletion, with $a$ still in $Q$.
Then, after the deletion, $J\cap Q=\{a\}$, and $J$ contains a vertex of $X\setminus\{u\}$.
If $a=u$, then $J$ and $F$ have a common vertex other than $u$.
Lemma~\ref{res:route-3-10}(i) and the definition of $\mathcal S$ give $(J,a,b)=(F,u,v_0)$.
If $a\ne u$, then $a\notin F$, since $F\cap Q=\{u\}$ after the deletion.
Also, $u\notin J$, since $u\in Q$ and $J\cap Q=\{a\}$.
The first paragraph gives $J\cap F=\varnothing$, a contradiction.
Deleting vertices from $Q$ cannot make a previously uncounted triple satisfy the counting conditions.
The count therefore decreases by $t+1$, which is at most $k$ if $t<k$, and equals
$k+1$ if $t=k$, precisely when a new set $K_i$ is constructed.
The initial count is $b_0$ and the final count is zero.
Summing over all steps gives $b_0\le p-|Q|+m\le p-1+m$, where $Q$ is the set remaining at the end.

We now prove $2m\le\sigma$. We may assume $m>0$, since $\sigma\ge0$.
First, let $S_u,S_w$ be small maximal nontrivial critical sets with distinct roots
$u,w\in P$ and respective designated neighbours $y(u),y(w)$.
Let $K$ be a vertex set such that $\{u,w\}\subseteq K\subseteq S_u\cap S_w$.
By interchanging $u$ and $w$, if necessary, we may assume that $|S_w|\le|S_u|$.
Let $v_1,v_2$ be the two neighbours of the cubic vertex $u$ other than $y(u)$,
and put $F_1=S_{u,v_1}$ and $F_2=S_{u,v_2}$.
By Lemma~\ref{res:route-3-10}(i), $S_u\cap F_1=S_u\cap F_2=\{u\}$.
Since $w\in S_u$ and $w\ne u$, neither $F_1$ nor $F_2$ contains $w$.
As $u\in F_1\cap S_w$ and $|F_2|\ge2$, Lemma~\ref{res:route-3-11}(iii),
applied to $S_u,F_1,F_2$, gives
\[
 |F_1\cup S_w|\le |F_1|+|S_w|-1
 \le |F_1|+|S_u|-1\le q-|F_2|-1\le q-3.
\]
Lemma~\ref{res:route-3-12}, applied to $S_w$ and $F_1$, gives
$F_1\subseteq S_w\setminus\{w\}$. Similarly, $F_2\subseteq S_w\setminus\{w\}$.
In particular, $F_1$ and $F_2$ are small.

Let $h\in N(K)\setminus\{y(u),y(w)\}$.
For each $x\in\{u,w\}$ adjacent to $h$, let $z$ be the neighbour of $x$ distinct from $y(x)$ and $h$,
and put $F_x=S_{x,z}$.
By Lemma~\ref{res:route-3-10}(i), $F_x\cap S_x=\{x\}$.
Since $x\in K\subseteq S_x$, we have $F_x\cap K=\{x\}$.
Since $h\in N(F_x)\setminus\{z\}$, Lemma~\ref{comp:lem:root-deletion} implies that $h$ has a neighbour
in $F_x\setminus\{x\}$, and hence outside $K$.
If both $u$ and $w$ are adjacent to $h$, then $F_u$ is one of $F_1,F_2$,
so $F_u\subseteq S_w\setminus\{w\}$.
By Lemma~\ref{res:route-3-10}(i), $F_w\cap S_w=\{w\}$.
Thus $F_u\cap F_w=\varnothing$, and $h$ has two distinct neighbours outside $K$,
one in $F_u$ and the other in $F_w$.

Let $Y_i=\{y_i(x):x\in R_i\}$.
Let $h$ be a vertex of $G$, and let $l$ be the number of sets $K_i$ such that
$h\in N(K_i)\setminus Y_i$. For each such set, put $d_i(h)=|N(h)\cap K_i|$
and $r_i(h)=|N(h)\cap R_i|$. We show that
\begin{equation}\label{eq:route-3-5}
 \sum_i\bigl(d_i(h)+r_i(h)-3\bigr)\le d(h)-3,
\end{equation}
where the sum is taken over the indices of these $l$ sets.
If $l=0$, this follows from $d(h)\ge3$.
Suppose $l>0$. Since $K_i\cap K_j=R_i\cap R_j$ and the incidence graph is a forest,
two distinct sets among these $l$ sets have at most one vertex in common,
and any common vertex belongs to both $R_i$ and $R_j$.
If $K_i\cap K_j=\{v\}$, then $|K_i\cup K_j|\le2b_q-1\le q-3$.
Lemma~\ref{fd:lem:uncrossing} therefore gives $N(K_i)\cap N(K_j)=N(v)$.
As $h$ belongs to both neighbourhoods, $hv\in E(G)$.
For the sets $K_i$ containing $v$, the vertices $y_i(v)$ are distinct neighbours of $v$
and are different from $h$. Since $v$ is cubic, at most two of the $l$ sets contain $v$.
Consider the graph with one vertex for each of these sets, where two distinct vertices
are adjacent if and only if the corresponding sets have a vertex in common.
Each $K_i$ can meet another set only at one of the two vertices of $R_i$,
and each of these vertices belongs to at most one other set.
Hence this graph has maximum degree at most two.
A cycle in this graph would give a cycle in the incidence graph, a contradiction.
Thus each component is a path, possibly consisting of a single vertex.

Let $c$ be the number of these paths.
On each path, exactly two roots belong to only one of its sets $K_i$.
Denote these roots by $u,w$, and let $K^*$ be the union of the sets on the path.
For $x\in\{u,w\}$, let $S_x$ be a small maximal nontrivial critical set
containing the set $K_i$ at that end of the path, with root $x$ and designated neighbour $y_i(x)$.
Such a set exists by the construction of $K_i$.
Every other set $K_j$ on the path excludes $x$.
Starting at $K_i$, apply Lemma~\ref{res:route-3-12} successively to $S_x$
and the next set on the path. These sets intersect and both have size at most $b_q$.
Thus $K^*\subseteq S_u\cap S_w$.
The preceding argument now gives a neighbour of $h$ outside $K^*$ for each of $u,w$
adjacent to $h$. If both are adjacent to $h$, these neighbours lie in disjoint sets
$F_u,F_w$, at least one of which is small.

Set $a=2l-\sum_i r_i(h)$, and let $s$ be the number of neighbours of $h$ outside the union of these $l$ sets.
Since $|R_i|=2$ for every $i$ and every root belonging to two sets is adjacent to $h$,
$a$ is the number of roots that belong to exactly one of these sets and are not adjacent to $h$.
The graph defined above has $l$ vertices and $c$ components, each of which is a path.
It therefore has $l-c$ edges, each corresponding to a root belonging to two sets.
Each such root is counted twice in $\sum_i d_i(h)$, while every other neighbour of $h$ in the union is counted once.
Hence $\sum_i d_i(h)=d(h)-s+l-c$. Using $\sum_i r_i(h)=2l-a$, we obtain
\[
 \sum_i\bigl(d_i(h)+r_i(h)-3\bigr)=d(h)-s-a-c.
\]
It suffices to show $s+a+c\ge3$.
This holds if $c\ge3$, or if $c=2$ and $a\ge1$.
If $c=1$, exactly $2-a$ of the two roots $u,w$ are adjacent to $h$.
By the preceding paragraph, $h$ has at least $2-a$ distinct neighbours outside $K^*$.
Since $K^*$ is the union of all $l$ sets in this case, $s\ge2-a$, as required.

Finally, suppose $c=2$ and $a=0$, and assume, for a contradiction, that $s=0$.
Let $K^*,L^*$ be the unions of the sets on the two paths. These unions are disjoint.
Let $u,w$ be the two roots that each belong to exactly one set on the first path.
Since $a=0$, both $u$ and $w$ are adjacent to $h$.
By the preceding argument, the disjoint sets $F_u,F_w$ each contain a neighbour of $h$ outside $K^*$.
As $s=0$, both neighbours belong to $L^*$.
By interchanging $u$ and $w$, if necessary, we may assume that $F_u$ is small.
Select a set $K_j$ on the second path with $K_j\cap F_u\ne\varnothing$.
Since $u\in K^*$ and $K^*\cap L^*=\varnothing$, we have $u\notin K_j$.
Also, $|F_u\cup K_j|\le2b_q-1\le q-3$, so Lemma~\ref{res:route-3-12} gives $K_j\subseteq F_u$.
If a set $K_i$ on the second path is adjacent to a set already contained in $F_u$,
then $K_i\cap F_u\ne\varnothing$.
Since $u\notin K_i$ and $|K_i|\le b_q$, the same lemma gives $K_i\subseteq F_u$.
Repeating this argument along the second path gives $L^*\subseteq F_u$.
This contradicts $F_u\cap F_w=\varnothing$, since $F_w\cap L^*\ne\varnothing$.
Thus $s\ge1$, so $s+a+c\ge3$, proving \eqref{eq:route-3-5}.

For each $i$, $Y_i$ consists of two vertices, each of which has exactly one neighbour in $K_i$.
Since $K_i$ is contained in one part of the bipartition, counting the edges incident with vertices of $K_i$ gives
\[
 \sum_{h\in N(K_i)\setminus Y_i}d_i(h)=\sum_{x\in K_i}d(x)-2.
\]
The two roots in $R_i$ are cubic, and exactly two edges join $R_i$ to $Y_i$. Hence
\[
 \sum_{h\in N(K_i)\setminus Y_i}r_i(h)=6-2=4.
\]
Since $|N(K_i)|=|K_i|+2$ and $|Y_i|=2$, we have $|N(K_i)\setminus Y_i|=|K_i|$.
Consequently, for each $i$,
\[
 \sum_{h\in N(K_i)\setminus Y_i}\bigl(d_i(h)+r_i(h)-3\bigr)
 =2+\sum_{x\in K_i}(d(x)-3).
\]
Summing \eqref{eq:route-3-5} over all $h\in V(G)$ and changing the order of summation gives
\[
 2m+\sum_i\sum_{x\in K_i}(d(x)-3)\le\sigma.
\]
Since $d(x)\ge3$ for every vertex $x$, we obtain $2m\le\sigma$.
Thus $b_0\le p-1+m\le p-1+\frac{\sigma}{2}$, as required.
\end{proof}

Let $\mathcal A$ consist of the small nontrivial critical sets $S_{u,v}$ such that $S_{u,v}\cap P=\{u\}$.
Let $\ell$ denote the number of triples $(S,u,v)\in\mathcal S$ such that $S$ is large.
The following theorem gives a lower bound on $|\mathcal A|$.

\begin{theorem}\label{res:route-3-19}
Let $G$ be a minimal simple brace of order $2q$, where $q\ge6$.
Two distinct members of $\mathcal A$ are disjoint if their roots are different.
If they have the same root $u$, then $u$ is their only common vertex.
If $p>0$, then
$|\mathcal A|\ge2p-\ell+1-\frac\sigma2$.
\end{theorem}
\begin{proof}
Let $S,T$ be distinct members of $\mathcal A$, with roots $u,v$, respectively.
Suppose that $u\ne v$. By the definition of $\mathcal A$, $u\notin T$ and $v\notin S$.
If $S\cap T\ne\varnothing$, then $|S\cup T|\le2b_q-1\le q-3$.
Lemma~\ref{res:route-3-12} gives $T\subseteq S$, contrary to $v\in T\setminus S$.
Hence $S$ and $T$ are disjoint.
If $u=v$, then $S$ and $T$ have different designated neighbours, since only one critical set
was selected for each root and designated neighbour.
By Lemma~\ref{res:route-3-10}(i), $u$ is their only common vertex.

For each $S\in\mathcal A$, there is exactly one triple $(S,u,v)\in\mathcal S$.
The root must be the unique vertex $u$ of $S\cap P$.
If $(S,u,v)$ and $(S,u,w)$ are two such triples with $v\ne w$,
then Lemma~\ref{res:route-3-10}(i) gives $S=\{u\}$, contrary to $|S|\ge2$.

Assume $p>0$, and let $b_0$ denote the number of triples $(S,u,v)\in\mathcal S$
such that $S$ is small and $|S\cap P|\ge2$.
By Lemma~\ref{res:route-3-18}, $b_0\le p-1+\frac{\sigma}{2}$.
Every triple $(S,u,v)\in\mathcal S$ satisfies $u\in S\cap P$.
Of the $3p$ triples, $\ell$ have $S$ large.
Among the triples with $S$ small, $b_0$ satisfy $|S\cap P|\ge2$.
The remaining triples have $S$ small and $S\cap P=\{u\}$, so their number is $|\mathcal A|$.
Hence
\begin{equation}\label{eq:route-3-6}
 \begin{aligned}
 |\mathcal A|&=3p-\ell-b_0\\
 &\ge2p-\ell+1-\frac\sigma2.
 \end{aligned}
\end{equation}
\end{proof}

\section{Thin Edges and Critical Sets}\label{sec:short-critical}
In this section, we study thin edges incident with vertices of small critical sets or with neighbours of these vertices.


\subsection{Thin Edges and Small Critical Sets}

We begin with some notation. Let $G$ be a minimal simple brace with bipartition $(A,B)$.
Recall that $C$ is the set of cubic vertices of $G$ incident with at least one thin edge.
A cubic vertex $x$ incident with no thin edge belongs to $D$ if and only if there exist three pairwise disjoint nontrivial critical sets, one rooted at each neighbour of $x$ and all with designated neighbour $x$.
The set $P$ consists of the cubic vertices outside $C\cup D$.
Put $V_{\ge4}=\{v\in V(G):d(v)\ge4\}$ and $U=V_{\ge4}\cup D$.
We need the following two lemmas.
\begin{lemma}\label{res:route-3-20}
Let $G$ be a minimal simple brace of order $2q$, where $q\ge6$.
If $x,h$ are distinct vertices of $G$ in the same part of its bipartition and
$N(x)\subseteq N(h)$, every edge joining $h$ to a vertex of $N(x)$ is thin.
\end{lemma}
\begin{proof}
Suppose $ha$ is nonthin for some $a\in N(x)$. By Proposition~\ref{res:route-3-1} and Lemma~\ref{res:route-3-2}(i),
there is a nontrivial critical set $R$ with root $a$ and designated neighbour $h$.
Since $N(h)\cap R=N(x)\cap R=\{a\}$, neither $h$ nor $x$ has a neighbour in $R\setminus\{a\}$.
Hence $|N(R\setminus\{a\})|\le|N(R)|-2=|R|$. This contradicts Theorem~\ref{thm:brace-hall}, since
$1\le|R\setminus\{a\}|\le q-4$.
\end{proof}

\begin{lemma}\label{lem:cubic-common-neighbours}
If two distinct cubic vertices of a minimal simple brace $G$, where $|V(G)|=2q$ and $q\ge6$,
have two common neighbours, then neither belongs to $D$.
\end{lemma}
\begin{proof}
Let $x,z$ be distinct cubic vertices with $N(x)=\{a,b,y\}$ and $N(z)=\{a,b,y'\}$.
By Theorem~\ref{thm:brace-hall}, $|N(\{x,z\})|\ge4$, so $y\ne y'$.
Thus $\{x,z\}$ is a nontrivial critical set with root $x$ and designated neighbour $y$.
By Lemma~\ref{res:route-3-2}(ii), let $R$ be the nontrivial critical set paired with $\{x,z\}$,
with root $y$ and designated neighbour $x$. Then $|R|=q-3$.
If $x\in D$, select nontrivial critical sets $R_a,R_b$ with respective roots $a,b$ and common designated neighbour $x$.
By Lemma~\ref{res:route-3-6}(ii) and the definition of $D$, the sets $R,R_a,R_b$ are pairwise disjoint.
All three sets are contained in the same part of the bipartition, which has $q$ vertices,
but their total size is at least $q-3+2+2>q$, a contradiction.
Thus $x\notin D$. By symmetry, $z\notin D$.
\end{proof}

A nontrivial critical set $S$ is \textit{small} if $|S|\le b_q$ and large otherwise, where $b_q=\lfloor\frac{q-2}{2}\rfloor$.
For each $u\in P$ and $v\in N(u)$, let $S_{u,v}$ denote the maximal nontrivial critical set
selected in Subsection~3.3 with root $u$ and designated neighbour $v$.
We shall use Lemmas~\ref{res:route-3-20} and~\ref{lem:cubic-common-neighbours} in the proof of the following lemma.
\begin{lemma}\label{res:route-3-21}
If $\{y,z\}$ is a small maximal nontrivial critical set with root $y$ and designated neighbour $x$
in a minimal simple brace $G$, where $|V(G)|=2q$ and $q\ge6$, then $z\in C\cup D$
and the following statements hold.
\begin{enumerate}[label=\textup{(\roman*)},leftmargin=*]
\item If $y$ is cubic, then $y\notin D$ and $z\in C$.
\item If $y$ is not cubic, then $d(y)=4$, $N(z)=N(y)\setminus\{x\}\subseteq C$,
and every edge joining $y$ to a vertex of $N(z)$ is thin.
\end{enumerate}
\end{lemma}
\begin{proof}
Put $S=\{y,z\}$. Since $x\in N(y)\setminus N(z)$, $|N(S)|=4$ and $G$ has minimum
degree three or more, $z$ is cubic and $d(y)\le4$. Suppose that $z\in P$.
By Lemma~\ref{res:route-3-10}(iii), at least two of the three sets $S_{z,v}$, where $v\in N(z)$, are small.
By Lemma~\ref{res:route-3-10}(i), these two sets have only $z$ in common.
At least one of them, say $F$, does not contain $y$.
As $|S\cup F|=|F|+1\le b_q+1\le q-3$, Lemma~\ref{res:route-3-12}, applied to $S$ and $F$,
gives $F\subseteq\{z\}$, contrary to $|F|\ge2$. Hence $z\in C\cup D$.

To prove (i), suppose that $y$ is cubic. Then $y$ and $z$ have two common neighbours.
By Lemma~\ref{lem:cubic-common-neighbours}, $y,z\notin D$. Since $z\in C\cup D$, it follows that $z\in C$.
To prove (ii), suppose that $y$ is not cubic. Then $d(y)=4$ and $N(z)=N(y)\setminus\{x\}$.
By Lemma~\ref{res:route-3-20}, every edge joining $y$ to a vertex of $N(z)$ is thin.
Since $d(y)=4$, Lemma~\ref{res:route-3-8}(i) implies that every vertex of $N(z)$ is cubic.
Hence $N(z)\subseteq C$.
\end{proof}

We give the following property of small maximal nontrivial critical sets.

\begin{lemma}\label{res:route-3-23}
Let $S$ be a small maximal nontrivial critical set with root $u$ and designated neighbour $v$, and let
$x\in S\setminus\{u\}$.
If $xy$ is a nonthin edge and $S$ contains no nontrivial critical set with root $x$ and designated neighbour $y$,
then there exists a small maximal nontrivial critical set $Q$ with root $y$ and designated neighbour $x$ such that
$Q\subseteq N(S)$ and $|N(Q)\setminus S|\le2$.
\end{lemma}
\begin{proof}
Interchanging $A$ and $B$ if necessary, assume that $S\subseteq A$.
Since $xy$ is nonthin, Proposition~\ref{res:route-3-1} and Lemma~\ref{res:route-3-2}(i)
give a nontrivial critical set with root $x$ and designated neighbour $y$. Select an inclusion-minimal such set $F$.
Suppose that $|S\cup F|\le q-3$. If $u\notin F$, Lemma~\ref{res:route-3-12} gives
$F\subseteq S\setminus\{u\}$, contrary to the hypothesis.
If $u\in F$, Lemma~\ref{fd:lem:uncrossing} gives $|N(S\cap F)|=|S\cap F|+2$.
The set $S\cap F$ contains the distinct vertices $u,x$, has size at most $q-3$ and
satisfies $N(y)\cap(S\cap F)=\{x\}$. Thus it is a nontrivial critical set with root $x$
and designated neighbour $y$ contained in $S$, again contrary to the hypothesis.
Hence $|S\cup F|\ge q-2$.

By Lemma~\ref{res:route-3-2}(ii), $F$ is paired with a nontrivial critical set $Q$
with root $y$ and designated neighbour $x$. By Lemma~\ref{res:route-3-2}(i),
$Q=(B\setminus N(F))\cup\{y\}$, $N(Q)=(A\setminus F)\cup\{x\}$ and $E(F,Q)=\{xy\}$.
Since $x\in S$, $N(Q)\setminus S=A\setminus(S\cup F)$ has at most two vertices.
Suppose that $z\in Q\setminus N(S)$. Then $z\ne y$, since $y\in N(S)$, and
$E(F,Q)=\{xy\}$ implies $N(z)\subseteq A\setminus(S\cup F)$.
Thus $d(z)\le2$, contrary to the fact that every vertex of $G$ has degree at least three.
Hence $Q\subseteq N(S)$.

Since $|S\cup F|\ge q-2$ and $x\in S\cap F$, we have
$|Q|=q-1-|F|\le|S|+1-|F\cap S|\le|S|$, so $Q$ is small.
Suppose that $Q$ is properly contained in a nontrivial critical set $Q'$ with root $y$
and designated neighbour $x$. By Lemma~\ref{res:route-3-2}(ii), (iii), the set $F'$ paired with $Q'$
has root $x$ and designated neighbour $y$ and satisfies $F'\subsetneq F$,
contrary to the choice of $F$. Thus $Q$ is maximal.
\end{proof}

The following lemma is a consequence of Lemma~\ref{res:route-3-13}.
\begin{lemma}\label{res:route-3-24}
Let $S$ be a small maximal nontrivial critical set with root $u$ and designated neighbour $v$, and let
$Q\subseteq N(S)\setminus\{v\}$ be a nontrivial critical set.
If $cz$ is a thin edge of $G$, where $c$ is cubic and $z\in Q$, then $c\in S$.
\end{lemma}
\begin{proof}
Suppose that $c\notin S$. If $N(c)\cap Q=\{z\}$, then $Q$ is a nontrivial critical set
with root $z$ and designated neighbour $c$. Lemma~\ref{res:route-3-2}(ii) and Proposition~\ref{res:route-3-1}
imply that $cz$ is nonthin, contrary to the hypothesis. Thus $c$ has at least two neighbours in $Q$.
By Lemma~\ref{res:route-3-13}, $c$ has exactly two neighbours in $N(S)$, one of which is $v$.
Since $Q\subseteq N(S)\setminus\{v\}$, $c$ has at most one neighbour in $Q$, a contradiction.
\end{proof}

We use Lemma~\ref{res:route-3-23} to prove the result below.

\begin{lemma}\label{res:route-3-26}
Let $S$ be a small maximal nontrivial critical set with root $u$ and designated neighbour $v$, and let
$x\in U\cap(S\setminus\{u\})$ be incident with at least three nonthin edges.
Consider all nonthin edges $xy$ for which no nontrivial critical set with root $x$ and designated neighbour $y$
is contained in $S$. The following statements hold.
\begin{enumerate}[label=\textup{(\roman*)},leftmargin=*]
\item For each such edge $xy$, there is a small maximal nontrivial critical set $Q_y$ with root $y$
and designated neighbour $x$ such that $Q_y\subseteq N(S)$, $|N(Q_y)\setminus S|\le2$
and $|Q_y|\le|S|-2$.
\item The sets $Q_y$ in \textup{(i)} can be selected so that they are pairwise disjoint and at least two of them do not contain $v$.
\end{enumerate}
\end{lemma}
\begin{proof}
Assume that $S\subseteq A$.
Every nontrivial critical set with root $x$ has at least three vertices.
Suppose that $\{x,w\}$ is such a set, with designated neighbour $y$.
Since $y\in N(x)\setminus N(w)$ and $|N(\{x,w\})|=4$, $w$ is cubic and $d(x)\le4$.
If $x$ is cubic, then $x\in D$, but $x,w$ have two common neighbours, contrary to Lemma~\ref{lem:cubic-common-neighbours}.
If $d(x)=4$, then $N(w)=N(x)\setminus\{y\}$.
By Lemma~\ref{res:route-3-20}, the three edges joining $x$ to $N(w)$ are thin.
Thus $x$ is incident with at most one nonthin edge, contrary to the hypothesis.

We next show that there is at most one nonthin edge $xy$ for which $S$ contains a nontrivial critical set
with root $x$ and designated neighbour $y$.
Let $F\subseteq S$ be such a set. If $u\in F$, then $N(v)\cap F=\{u\}$, and
Lemma~\ref{comp:lem:root-deletion} shows that $F\setminus\{u\}$ is a nontrivial critical set with root $x$
and designated neighbour $y$, since $|F|\ge3$ and $x\ne u$.
We may therefore assume that $F\subseteq S\setminus\{u\}$.
By Lemma~\ref{res:route-3-2}(ii), its paired set $R_y\subseteq B$ has root $y$
and designated neighbour $x$, with $|R_y|=q-1-|F|\ge q-|S|$.
For two such edges, Lemma~\ref{res:route-3-8}(ii) would give disjoint sets $R_{y_1},R_{y_2}\subseteq B$
of total size at least $2(q-|S|)>q$, a contradiction.

For each nonthin edge $xy$ considered in (i), apply Lemma~\ref{res:route-3-23} to select a small maximal nontrivial critical set $Q_y$
with root $y$ and designated neighbour $x$ such that $Q_y\subseteq N(S)$ and $|N(Q_y)\setminus S|\le2$.
By Lemma~\ref{res:route-3-8}(ii), the sets $Q_y$ are pairwise disjoint.
If $S$ contains no nontrivial critical set with root $x$, there are at least three sets $Q_y$.
For each $Q_y$, two other sets contain at least four vertices of $N(S)\setminus Q_y$.
Hence $|Q_y|\le|N(S)|-4=|S|-2$.
At most one $Q_y$ contains $v$, so at least two do not contain $v$.

Otherwise, let $F\subseteq S\setminus\{u\}$ be a nontrivial critical set with root $x$ and designated neighbour $y_0$. Its paired set is $R=(B\setminus N(F))\cup\{y_0\}$.
Since $y\ne y_0$, Lemma~\ref{res:route-3-8}(ii) gives $Q_y\cap R=\varnothing$.
Thus $Q_y\subseteq B\setminus R=N(F)\setminus\{y_0\}$, a set of $|F|+1$ vertices.
As at most one nonthin edge $xy$ is excluded, there are at least two sets $Q_y$.
Each has at least two vertices, and they are pairwise disjoint.
Hence $|Q_y|\le|F|-1\le|S|-2$.
Finally, $v\notin N(F)$, since $N(v)\cap S=\{u\}$ and $u\notin F$.
Hence no $Q_y$ contains $v$. This proves (i) and (ii).
\end{proof}

For $R\subseteq E(G)$ and $v\in V(G)$, let $N_R(v)$ denote the set of vertices joined to $v$ by an edge of $R$.
For $X\subseteq V(G)$, let $N_R(X)$ denote the set of vertices outside $X$ joined to a vertex of $X$ by an edge of $R$.
Put
$C_0=\{c\in C:N_{T(G)}(c)\cap C=\varnothing\}$ and $C_+=C\setminus C_0$.
A vertex $p_0\in P$ is an \emph{owner} of $c\in C$ if, for some $v\in N(p_0)$, $\{p_0,c\}$ is a maximal nontrivial
critical set with root $p_0$ and designated neighbour $v$.
In this case, we say that $c$ is \emph{owned} by $p_0$.

\begin{lemma}\label{res:route-3-28}
Let $G$ be a minimal simple brace, where $|V(G)|=2q$ and $q\ge6$, and let $p_0\in P$ and $c\in C$ be adjacent. For every thin edge $cz$, the vertex $z$
belongs to $C$ and is owned by $p_0$.
Moreover, $E(P,C_0)=\varnothing$ and $|N(c)\cap P|\le1$.
\end{lemma}
\begin{proof}
Let $p_0c\in E(P,C)$, and write $N(p_0)=\{c,a,b\}$ and $N(c)=\{p_0,r,s\}$.
Since $p_0\in P$, the edges $p_0a,p_0b$ are nonthin. By Proposition~\ref{res:route-3-1} and Lemma~\ref{res:route-3-2}(i),
select maximal nontrivial critical sets $S_a,S_b$ with root $p_0$ and respective designated neighbours $a,b$.
Each of $S_a,S_b$ contains a vertex of $\{r,s\}$. If, say, $S_a\cap\{r,s\}=\varnothing$, then
$N(S_a\setminus\{p_0\})\subseteq N(S_a)\setminus\{a,c\}$, so $|N(S_a\setminus\{p_0\})|\le|S_a|$,
contrary to Theorem~\ref{thm:brace-hall}.
By Lemma~\ref{res:route-3-10}(i), $S_a\cap S_b=\{p_0\}$. Adjust notation so that $r\in S_a$ and $s\in S_b$.

Since $p_0c$ is nonthin, every thin edge incident with $c$ is $cr$ or $cs$.
By symmetry, suppose that $cr$ is thin. If $|S_a|\ge3$, Lemma~\ref{comp:lem:root-deletion} gives
$|N(S_a\setminus\{p_0\})|=|S_a|+1$, and $2\le|S_a\setminus\{p_0\}|\le q-5$ by Lemma~\ref{res:route-3-11}(iii).
Since $N(c)\cap(S_a\setminus\{p_0\})=\{r\}$, the set $S_a\setminus\{p_0\}$ is a nontrivial critical set
with root $r$ and designated neighbour $c$. By Proposition~\ref{res:route-3-1} and Lemma~\ref{res:route-3-2}(ii),
$cr$ is nonthin, a contradiction.
Thus $S_a=\{p_0,r\}$ is a small maximal nontrivial critical set, and Lemma~\ref{res:route-3-21}(i) gives $r\in C$.
By the definition of owner, $r$ is owned by $p_0$.
The same argument applies to $cs$ if it is thin. Since $c\in C$, it follows that $c\notin C_0$.
As $p_0c$ was arbitrary, $E(P,C_0)=\varnothing$. Finally, Lemma~\ref{res:route-3-17} gives $|N(c)\cap P|\le1$.
\end{proof}

We next consider critical sets containing vertices that have an owner.
Here, $\mathcal A$ consists of the small nontrivial critical sets $S_{u,v}$ satisfying $S_{u,v}\cap P=\{u\}$.
\begin{lemma}\label{res:route-3-29}
Let $\{p_0,c\}$ be a maximal nontrivial critical set with root $p_0\in P$ and designated neighbour $v$, where $c\in C$.
Then $p_0$ is the unique owner of $c$,
and $\{p_0,c\}$ is the only nontrivial critical set with root $p_0$ and designated neighbour $v$.
Every small maximal nontrivial critical set containing $c$ contains $p_0$,
and every member of $\mathcal A$ containing $c$ equals $\{p_0,c\}$.
\end{lemma}
\begin{proof}
Put $K=\{p_0,c\}$, and let $F$ be another
nontrivial critical set with root $p_0$ and designated neighbour $v$. Lemma~\ref{res:route-3-11}(iii) gives $|F|\le q-4$, so
$|K\cup F|\le q-3$. By Lemma~\ref{fd:lem:uncrossing}, $|N(K\cup F)|=|K\cup F|+2$.
Since $N(v)\cap(K\cup F)=\{p_0\}$, it is a nontrivial critical set with root $p_0$ and designated neighbour $v$.
The maximality of $K$ gives $F\subseteq K$, and $|F|\ge2$ gives $F=K$.

If a small maximal set $S$ contains $c$ but not $p_0$, then
$|K\cup S|\le b_q+1\le q-3$, and Lemma~\ref{res:route-3-12} gives $S\subseteq\{c\}$, a contradiction.
Thus every such $S$ contains $p_0$. If $p_1$ is an owner of $c$, applying this to $S=\{p_1,c\}$ gives $p_1=p_0$.
If $S\in\mathcal A$, it is rooted at $p_0$. Since $c\in S\cap K$,
Lemma~\ref{res:route-3-10}(i) implies that $v$ is also the designated neighbour of $S$.
Thus $S=K$ by the first paragraph.
\end{proof}

\Needspace{6\baselineskip}
We next consider the neighbours of a vertex and its owner.
\begin{lemma}\label{route:lem:owner-pair}
Let $p_0$ be an owner of $c$, and write $N(p_0)=\{a,b,v\}$ and $N(c)=\{a,b,w\}$.
The following statements hold.
\begin{enumerate}[label=\textup{(\roman*)},leftmargin=*]
\item The edge $cw$ is nonthin.
\item If $cb$ is nonthin, then $b\in P\cup C$, and there are
vertices $z,h$ such that $N(b)=\{p_0,c,z\}$ and $N(v)=\{p_0,z,h\}$.
The edges $ca$ and $vz$ are the unique thin edges incident with $c$ and $v$, respectively.
Moreover, $\{b,v\}$ is a small maximal nontrivial critical set with root $b$ and designated neighbour $c$.
\item If $cb$ is nonthin and $c\in C_+$, then $b\in P$ and $a$ is owned by $b$.
\end{enumerate}
\end{lemma}
\begin{proof}
The set $\{p_0,c\}$ is a nontrivial critical set with root $c$ and designated neighbour $w$, so $cw$ is
nonthin by Lemma~\ref{res:route-3-2}(ii) and Proposition~\ref{res:route-3-1}. Suppose that $cb$ is nonthin.
Since $cb$ and $cw$ are nonthin and $c\in C$, the edge $ca$ is the unique thin edge incident with $c$.
Since $bp_0\in E(G)$, no subset of $\{p_0,c\}$ is a nontrivial critical set with root $c$ and designated neighbour $b$.
By Lemma~\ref{res:route-3-23}, there is a small maximal nontrivial critical set $Q$ with root $b$ and designated
neighbour $c$ such that $Q\subseteq N(\{p_0,c\})$.
Since $N(c)\cap Q=\{b\}$, we have $Q\subseteq\{b,v\}$. As $|Q|\ge2$, it follows that $Q=\{b,v\}$.
By Lemma~\ref{res:route-3-21}, $v$ is cubic. If $b$ were not cubic, Lemma~\ref{res:route-3-21}(ii) would give
$N(v)=N(b)\setminus\{c\}$ and imply that $bp_0$ is thin, contrary to $p_0\in P$.
Thus $b$ is cubic. By Lemma~\ref{res:route-3-21}(i), $b\notin D$ and $v\in C$, so $b\in P\cup C$.

Write $N(b)=\{p_0,c,z\}$. Since $|N(\{b,v\})|=4$, the cubic vertices $b,v$ have exactly two common neighbours.
As $c\notin N(v)$, these common neighbours are $p_0,z$. Thus $N(v)=\{p_0,z,h\}$ with $h\notin\{p_0,c,z\}$.
The edge $vp_0$ is nonthin because $p_0\in P$. The set $\{b,v\}$ is a nontrivial critical set with root $v$
and designated neighbour $h$, so $vh$ is nonthin by Lemma~\ref{res:route-3-2}(ii) and Proposition~\ref{res:route-3-1}.
As $v\in C$, the edge $vz$ is the unique thin edge incident with $v$.

Suppose that $c\in C_+$. Since $ca$ is the unique thin edge incident with $c$, the definition of $C_+$ gives $a\in C$.
Write $N(a)=\{p_0,c,a'\}$. Suppose that $b\in C$. Since $bp_0$ and $bc$ are nonthin, $bz$ is thin.
By Theorem~\ref{thm:brace-hall}, $|N(\{a,b\})|\ge4$, so $a'\ne z$.
Thus $\{a,b\}$ is a nontrivial critical set with root $b$ and designated neighbour $z$.
By Lemma~\ref{res:route-3-2}(ii) and Proposition~\ref{res:route-3-1}, $bz$ is nonthin, a contradiction.
Hence $b\in P$. Since $ca$ is thin, Lemma~\ref{res:route-3-28} applied to $bc$ implies that $b$ is an owner of $a$.
\end{proof}

We introduce one further definition.
A \emph{$2$-bounded set} is a set $X$ contained in one part of the bipartition such that $|X|\ge2$, $|N(X)|=|X|+2$,
and every vertex of $X$ is incident with at most two nonthin edges.
The next two lemmas will be used in the counting arguments.
\begin{lemma}\label{res:route-3-32}
Let $G$ be a minimal simple brace of order $2q$, where $q\ge6$, and let $S$ be a small maximal nontrivial critical set of $G$ containing no $2$-bounded set.
If no vertex of $S\cap C$ has an owner, then there are two disjoint small maximal nontrivial
critical sets $Q_1,Q_2$ such that $Q_i\subseteq N(S)$, $|N(Q_i)\setminus S|\le2$
and $|Q_i|\le|S|-2$ for $i=1,2$.
\end{lemma}
\begin{proof}
Select a small maximal nontrivial critical set $R\subseteq S$ of minimum cardinality, and let $r$ be its root.
Suppose that $p_0\in(R\setminus\{r\})\cap P$. By Lemma~\ref{res:route-3-10}(iii), there are two small maximal
nontrivial critical sets with root $p_0$ and distinct designated neighbours. By Lemma~\ref{res:route-3-10}(i),
their intersection is $\{p_0\}$, so at least one of them, say $F$, does not contain $r$.
Since $p_0\in R\cap F$ and $|R\cup F|\le2b_q-1\le q-3$, Lemma~\ref{res:route-3-12} gives
$F\subseteq R\setminus\{r\}$. Thus $|F|<|R|$, contrary to the choice of $R$.
Hence no vertex of $R\setminus\{r\}$ belongs to $P$.

Suppose that $R=\{r,z\}$. By Lemma~\ref{res:route-3-21}, $z\in C\cup D$.
If $r$ is cubic, Lemma~\ref{res:route-3-21}(i) gives $r\in P\cup C$ and $z\in C$.
If $r\in P$, then $r$ is an owner of $z$, contrary to the hypothesis.
If $r\in C$, then $R$ is a $2$-bounded set contained in $S$, again a contradiction.
Suppose instead that $r$ is not cubic. By Lemma~\ref{res:route-3-21}(ii), $d(r)=4$ and three edges incident with $r$ are thin.
If $z\in C$, then $R$ is again a $2$-bounded set contained in $S$, a contradiction.
Thus $z\in D\subseteq U$, and all three edges incident with $z$ are nonthin.
By Lemma~\ref{res:route-3-26}(i) and (ii), applied to $R$ and $z$, there is a nontrivial critical set of size at most
$|R|-2=0$, whereas every nontrivial critical set has at least two vertices. Hence $|R|\ge3$.

By Lemma~\ref{comp:lem:root-deletion}, $|N(R\setminus\{r\})|=|R\setminus\{r\}|+2$.
Since $|R|\ge3$, the set $R\setminus\{r\}$ has at least two vertices. It is not a $2$-bounded set, since it is contained in $S$.
Hence some $x\in R\setminus\{r\}$ is incident with at least three nonthin edges.
Vertices of $C$ are incident with at most two nonthin edges, so $x\notin C$. Since $x\notin P$ by the first paragraph, we have $x\in U$.
By Lemma~\ref{res:route-3-26}(i) and (ii), there are two disjoint small maximal nontrivial critical sets
$Q_1,Q_2\subseteq N(R)$ such that $|N(Q_i)\setminus R|\le2$ and $|Q_i|\le|R|-2$ for $i=1,2$.
Since $R\subseteq S$, the sets $Q_1,Q_2$ are contained in $N(S)$, and each has at most two neighbours outside $S$
and at most $|S|-2$ vertices.
\end{proof}

\begin{lemma}\label{res:route-3-33}
Let $G$ be a minimal simple brace of order $2q$, where $q\ge6$, and let $S$ and $Q$ be small maximal nontrivial critical sets of $G$ such that $Q\subseteq N(S)$,
$|N(Q)\setminus S|\le2$ and $|Q|\le|S|-2$.
If no vertex of $(S\cup Q)\cap C$ has an owner, then at least one of the sets $S$ and $Q$ contains a $2$-bounded set.
\end{lemma}
\begin{proof}
Select a counterexample $(S,Q)$ for which $S$ has minimum cardinality, and let $u$ be the root of $S$.
Neither $S$ nor $Q$ contains a $2$-bounded set. By Lemma~\ref{res:route-3-32}, applied to $Q$, there are two disjoint
small maximal nontrivial critical sets $R_1,R_2\subseteq N(Q)$ such that $|N(R_i)\setminus Q|\le2$ and
$|R_i|\le|Q|-2$ for $i=1,2$.
We show that at least one of $R_1,R_2$ is contained in $S$.
If one of them, say $R_1$, is disjoint from $S$, then $R_1\subseteq N(Q)\setminus S$.
Since $|R_1|\ge2$ and $|N(Q)\setminus S|\le2$, we have $R_1=N(Q)\setminus S$.
As $R_1$ and $R_2$ are disjoint, it follows that $R_2\subseteq S$.
Otherwise, both $R_1$ and $R_2$ meet $S$. Since $R_1$ and $R_2$ are disjoint, select $i\in\{1,2\}$ such that $u\notin R_i$.
Then $|S\cup R_i|\le2b_q-1\le q-3$, so Lemma~\ref{res:route-3-12} gives $R_i\subseteq S\setminus\{u\}$.
In either case, select $i\in\{1,2\}$ with $R_i\subseteq S$.
Since $R_i\subseteq S$, it contains no $2$-bounded set, and no vertex of $(Q\cup R_i)\cap C$ has an owner.
Thus $(Q,R_i)$ satisfies the hypotheses of the lemma, but neither $Q$ nor $R_i$ contains a $2$-bounded set.
It is therefore another counterexample, with $|Q|\le|S|-2<|S|$, contrary to the choice of $S$.
\end{proof}

\subsection{Cubic Vertices and Large Critical Sets}

A \emph{domino} is a graph consisting of two $4$-cycles whose intersection consists of one edge and its endvertices.
We call the end vertices of their common edge the \emph{centres} and the remaining four vertices the \emph{corners} of the domino.

\Needspace{5\baselineskip}
\begin{lemma}\label{res:route-3-34}
Let $G$ be a minimal simple brace of order $2q$, where $q\ge6$. The following statements hold.
\begin{enumerate}[label=\textup{(\roman*)},leftmargin=*]
\item For each edge $uv\in E(G[P])$, there is a subgraph $H$ of $G$ that is a domino with centres $u,v$ and all four corners in $C$.
Each edge of $H$ joining two corners is the unique thin edge of $G$ incident with each of its end vertices.
\item The subgraphs in \textup{(i)} corresponding to distinct edges of $G[P]$ are vertex-disjoint.
Every vertex $c\in C_+$ with $d_{T(G)}(c)=1$ that has an owner is a corner of one of these subgraphs.
\end{enumerate}
\end{lemma}
\begin{proof}
We first prove (i). Let $uv\in E(G[P])$, with $u\in A$ and $v\in B$, and write
$N(u)=\{v,a,b\}$ and $N(v)=\{u,c,d\}$. Since $u,v\in P$, every edge incident with $u$ or $v$ is nonthin.
By Proposition~\ref{res:route-3-1} and Lemma~\ref{res:route-3-2}, select paired nontrivial critical sets $S_0,R_0$
with respective roots $u,v$ and respective designated neighbours $v,u$.
Also select maximal nontrivial critical sets $S_a,S_b$ with root $u$ and respective designated neighbours $a,b$,
and $R_c,R_d$ with root $v$ and respective designated neighbours $c,d$.

Each of $S_a,S_b$ contains a vertex of $\{c,d\}$. Otherwise, we may assume that $S_a\cap\{c,d\}=\varnothing$.
Then $N(a)\cap S_a=\{u\}$ and $N(v)\cap S_a=\{u\}$, so
$|N(S_a\setminus\{u\})|\le|N(S_a)|-2=|S_a|$.
This contradicts Theorem~\ref{thm:brace-hall}, since $1\le|S_a\setminus\{u\}|\le q-4$.
By Lemma~\ref{res:route-3-10}(i), $S_a\cap S_b=\{u\}$.
Adjust notation so that $c\in S_a$ and $d\in S_b$.
A similar argument shows that each of $R_c,R_d$ contains exactly one of $a,b$, and these vertices are distinct.

By Lemma~\ref{res:route-3-11}(ii), the sets $N(R_c)\setminus N(v)$ and $N(R_d)\setminus N(v)$
are disjoint from $N(R_0)\setminus N(v)$.
By Lemma~\ref{res:route-3-2}(i), $N(R_0)=(A\setminus S_0)\cup\{u\}$, so
$(N(R_c)\cup N(R_d))\setminus N(v)\subseteq S_0\setminus\{u\}$.
Since $\{a,b\}\subseteq R_c\cup R_d$, it follows that $N(a)\cup N(b)\subseteq N(v)\cup S_0$.
Similarly, $N(c)\cup N(d)\subseteq N(u)\cup R_0$.

Put $X=S_a\setminus\{u,c\}$. Since $N(a)\cap S_a=\{u\}$ and $N(v)\cap S_a=\{u,c\}$,
no vertex of $X$ is adjacent to $a$ or $v$.
By Lemma~\ref{res:route-3-10}(i), $S_a\cap S_0=\{u\}$.
Together with $N(b)\subseteq N(v)\cup S_0$, this gives $N(b)\cap S_a\subseteq\{u,c\}$,
so no vertex of $X$ is adjacent to $b$.
Thus $N(X)\subseteq N(S_a)\setminus N(u)$ and
$|N(X)|\le|S_a|+2-3=|X|+1$.
If $X\ne\varnothing$, then $1\le|X|\le q-5$, so Theorem~\ref{thm:brace-hall} gives $|N(X)|\ge|X|+2$, a contradiction.
Hence $S_a=\{u,c\}$, and similarly $S_b=\{u,d\}$.
The same argument at $v$ shows that $R_c,R_d$ also have two vertices.
These four maximal nontrivial critical sets are small because $q\ge6$.
Since their roots $u,v$ are cubic, Lemma~\ref{res:route-3-21}(i) gives $a,b,c,d\in C$.

Since $N(a)\cap S_a=\{u\}$ and $N(b)\cap S_b=\{u\}$, we have $ac,bd\notin E(G)$.
As $|N(\{u,c\})|=4$ and $|N(u)|=3$, the vertex $c$ has exactly one neighbour $w$ outside $N(u)$.
Since $c$ is cubic and is adjacent to $v$ but not to $a$, it follows that $N(c)=\{v,b,w\}$.
Similarly, $ad\in E(G)$.
The edge $cv$ is nonthin because $v\in P$.
Since $w\notin N(u)$, the set $\{u,c\}$ is a nontrivial critical set with root $c$ and designated neighbour $w$.
By Lemma~\ref{res:route-3-2}(ii) and Proposition~\ref{res:route-3-1}, $cw$ is nonthin.
As $c\in C$, the edge $bc$ is the unique thin edge incident with $c$.
A similar argument at $a,b,d$ shows that $bc,ad$ are the unique thin edges incident with their respective end vertices.
The cycles $uvcbu$ and $uvdau$ form a subgraph $H$ of $G$ that is a domino with centres $u,v$. This proves (i).

To prove (ii), consider two subgraphs given by (i) for distinct edges of $G[P]$.
Each centre is cubic and has exactly one neighbour in $P$, so the two subgraphs cannot have a common centre.
If they had a common corner $z$, then $z$ would have two distinct neighbours in $P$, contrary to Lemma~\ref{res:route-3-17}.
Since $P\cap C=\varnothing$, a centre of one subgraph cannot be a corner of the other.
Hence the two subgraphs are vertex-disjoint.

Finally, let $z\in C_+$ satisfy $d_{T(G)}(z)=1$, and let $p_0$ be an owner of $z$.
By Lemma~\ref{route:lem:owner-pair}(i), the edge joining $z$ to its neighbour outside $N(p_0)$ is nonthin.
Since $d_{T(G)}(z)=1$, there is a common neighbour $y$ of $p_0,z$ such that $zy$ is nonthin.
Lemma~\ref{route:lem:owner-pair}(iii) now gives $y\in P$.
By (i), the edge $p_0y$ corresponds to a subgraph of $G$ that is a domino with centres $p_0,y$.
As $y$ has degree three both in $G$ and in this subgraph, all its neighbours belong to the subgraph.
In particular, $z$ is one of its corners.
\end{proof}

The following lemma will be used in the proof of Lemma~\ref{res:route-3-37}.

\begin{lemma}\label{res:route-3-35}
Let $G$ be a minimal simple brace, and let $p_0\in P$ be an owner of $c\in C$.
Let $w$ be the unique neighbour of $c$ outside $N(p_0)$.
If $S\ne\{p_0,c\}$ is a nontrivial critical set with root $p_0$, then
$S\cup\{c\}$ is a nontrivial critical set with root $c$ and designated neighbour $w$.
\end{lemma}
\begin{proof}
Write $N(p_0)=\{a,b,v\}$ and $N(c)=\{a,b,w\}$, where $v$ is the designated neighbour
of $\{p_0,c\}$ with root $p_0$. Since $S\ne\{p_0,c\}$, Lemma~\ref{res:route-3-29} implies that the designated neighbour of $S$ is not $v$.
By Lemma~\ref{res:route-3-10}(i), $c\notin S$.
By Lemma~\ref{res:route-3-10}(ii), the sets $N(S)\setminus N(p_0)$ and
$N(\{p_0,c\})\setminus N(p_0)=\{w\}$ are disjoint.
Since $w\notin N(p_0)$, it follows that $w\notin N(S)$.
Also, $p_0\in S$, so both $a$ and $b$ belong to $N(S)$.
Hence $N(S\cup\{c\})=N(S)\cup\{w\}$ and $|N(S\cup\{c\})|=|N(S)|+1=|S\cup\{c\}|+2$.
Moreover, $N(w)\cap(S\cup\{c\})=\{c\}$, since $w\notin N(S)$ and $cw\in E(G)$.
Finally, Lemma~\ref{res:route-3-11}(iii) gives $|S|\le q-4$, so $2\le|S\cup\{c\}|\le q-3$.
Thus $S\cup\{c\}$ is a nontrivial critical set with root $c$ and designated neighbour $w$.
\end{proof}

We use Lemmas~\ref{res:route-3-34} and~\ref{res:route-3-35} to prove the following result.
\begin{lemma}\label{res:route-3-37}
Let $G$ be a minimal simple brace of order $2q$, where $q\ge6$.
For each large maximal nontrivial critical set $S=S_{p_0,b}$ with $p_0,b\in P$,
let $c_1,c_2$ be the two neighbours of $b$ other than $p_0$.
The following statements hold.
\begin{enumerate}[label=\textup{(\roman*)},leftmargin=*]
\item For $i=1,2$, the unique neighbour $w_i$ of $c_i$ outside $N(p_0)$ belongs to $U\cup C$.
The set $S\cup\{c_i\}$ is a large nontrivial critical set with root $c_i$ and designated neighbour $w_i$.
\item The vertices $c_1,c_2$ corresponding to all the large sets $S_{p_0,b}$ with $p_0,b\in P$ are pairwise distinct.
If $w_i\in C$, then $c_iw_i$ is nonthin. All pairs $(c_i,c_iw_i)$ with $w_i\in C$ are distinct.
\end{enumerate}
\end{lemma}
\begin{proof}
Let $S=S_{p_0,b}$ be large, where $p_0,b\in P$.
By Lemma~\ref{res:route-3-34}(i), $G$ has a subgraph $H$ that is a domino with centres $p_0,b$.
Since $b$ has degree three in both $G$ and $H$, the vertices $c_1,c_2$ are the two corners of $H$ adjacent to $b$.
Each $c_i$ belongs to $C$ and is joined by a thin edge to a vertex of $N(p_0)\cap C$.
Lemma~\ref{res:route-3-28} therefore implies that $p_0$ is an owner of each $c_i$.
Since $S$ is large, $S\ne\{p_0,c_i\}$.
By Lemma~\ref{res:route-3-35}, $S\cup\{c_i\}$ is a nontrivial critical set with root $c_i$ and designated neighbour $w_i$.
It is large because it contains $S$.
Also, $w_i\ne b$, since $w_i\notin N(p_0)$ and $b\in N(p_0)$.
As $c_i$ is cubic and has a neighbour $b\in P$, Lemma~\ref{res:route-3-17} gives $w_i\notin P$.
Hence $w_i\in U\cup C$. This proves (i).

By Lemma~\ref{res:route-3-11}(iii) and the definition of a large set, for each $p_0\in P$,
at most one of the three sets $S_{p_0,v}$, where $v\in N(p_0)$, is large.
Subgraphs corresponding to distinct edges of $G[P]$ are vertex-disjoint by Lemma~\ref{res:route-3-34}(ii).
If both $S_{p_0,b}$ and $S_{b,p_0}$ are large, the two vertices obtained from $S_{p_0,b}$ lie in the part of the bipartition containing $p_0$,
whereas those obtained from $S_{b,p_0}$ lie in the other part.
Thus all the vertices $c_i$ are distinct.
If $w_i\in C$, then $c_iw_i$ is nonthin by Lemma~\ref{route:lem:owner-pair}(i).
The pairs $(c_i,c_iw_i)$ with $w_i\in C$ are distinct because the vertices $c_i$ are distinct. This proves (ii).
\end{proof}

We conclude this section with a consequence of Lemma~\ref{res:route-3-9}.
\begin{lemma}\label{res:route-3-36}
Let $G$ be a minimal simple brace, and let $u_1v_1,u_2v_2\in E(G)$, where $u_1,u_2\in V(G)\setminus U$ are distinct
and $v_1,v_2\in U$. If for each $i=1,2$, there is a large nontrivial critical set
with root $u_i$ and designated neighbour $v_i$,
then $v_1$ and $v_2$ belong to distinct components of $G[U]$.
\end{lemma}
\begin{proof}
Suppose that $v_1,v_2$ lie in the same component of $G[U]$.
Since $u_1,u_2$ are distinct vertices outside $U$, there is a path from $u_1$ to $u_2$
starting with the edge $u_1v_1$ and ending with the edge $v_2u_2$ whose internal vertices belong to $U$.
Let $s$ be its length, so $s\ge2$.
Select large nontrivial critical sets $F_0,F_s$ with respective roots $u_1,u_2$ and respective designated neighbours $v_1,v_2$.
By Lemma~\ref{res:route-3-2}(ii) and Proposition~\ref{res:route-3-1}, the edges $u_1v_1,u_2v_2$ are nonthin.
By Lemma~\ref{res:route-3-8}(i), every thin edge has an end vertex in $C$.
Since $U\cap C=\varnothing$, every edge of $G[U]$ is nonthin.
Thus every edge of the path is nonthin.
The sets $F_0$ and $F_s$ satisfy
$|F_0|+|F_s|\ge2(b_q+1)\ge q-1$, so Lemma~\ref{res:route-3-9} gives
$q\ge|F_0|+|F_s|+s\ge q+1$, a contradiction.
\end{proof}


\section{The Number of Thin Edges in Braces}\label{sec:lower-bounds}


In this section, we prove our main results: Theorems~\ref{thm:main} and~\ref{thm:general}. We first establish the lower bound of thin edges for simple braces, and then extend our arguments to the case where parallel edges are allowed.

\subsection{The Proof of Theorem~\ref{thm:main}}\label{sec:proof-main}
By Lemma~\ref{prop:spanning-minimal-brace}, it suffices to prove Theorem~\ref{thm:main} for minimal simple braces.
Let $G$ be a minimal simple brace of order $n=2q$, where $q\ge6$.
We first consider members of $\mathcal A$ with at least three vertices.

Let $\mathcal A^*=\{S\in\mathcal A:|S|\ge3\}$.
For $c\in C$, let $\mu(c)=d_{T(G)}(c)-1$ and define
\[
\varepsilon(c)=
\begin{cases}
1, & \text{if } d_{T(G)}(c)=3,\\
0, & \text{otherwise}.
\end{cases}
\]
For $X\subseteq C$ contained in one part of the bipartition, let
$\alpha(X)=\sum_{z\in C\cap(X\cup N_{T(G)}(X))}\mu(z)$.
For $S\in\mathcal A^*$, set $w_S=|S\cap C_0|+\frac12|N(S)\cap C_0|$.
Let $f_S$ denote the number of nonthin edges with both end vertices in $C$ and one end vertex in $S$,
and let $e_S$ denote the number of thin edges joining $S\cap V_{\ge4}$ to $C_+$.

\begin{lemma}\label{res:route-4-1}
For every $S\in\mathcal A^*$,
\[
2\le w_S+f_S+2\left(\alpha(S\cap C)-\sum_{c\in S\cap C}\varepsilon(c)\right)+2\sum_{c\in N(S)\cap C_0}\mu(c)+2e_S.
\]
\end{lemma}
\begin{proof}
By the definitions of $\alpha$, $\mu$ and $\varepsilon$, $\alpha(S\cap C)\ge2\sum_{c\in S\cap C}\varepsilon(c)$.
Thus the three terms with coefficient two are nonnegative even integers.
If any of these terms is positive, the inequality follows. We may thus assume that all three are zero.
In particular, $\alpha(S\cap C)=\sum_{c\in S\cap C}\varepsilon(c)$, so $\alpha(S\cap C)=0$.
Also, $e_S=0$, and every vertex of $N(S)\cap C_0$ is incident with exactly one thin edge.
Since $\alpha(S\cap C)=0$, all vertices of $S\cap C$ and all vertices of $C$ joined to them by thin edges
are incident with exactly one thin edge.
In particular, the thin edges incident with vertices of $S\cap C_+$ have distinct end vertices outside $S$.
It remains to prove $w_S+f_S\ge2$ under these conditions.

Let $u$ be the root of $S$ and $v$ its designated neighbour.
By Lemma~\ref{res:route-3-29}, no vertex of $S\cap C$ has an owner, since $|S|\ge3$.
If a set $R$ contained in a part of the bipartition satisfies $|N(R)|=|R|+2$ and $2\le|R|\le q-3$, and $zc$ is thin
with $z\in R$, then $c$ has another neighbour in $R$; otherwise
$R$ would be a nontrivial critical set with root $z$ and designated neighbour $c$, contrary to Proposition~\ref{res:route-3-1}.

Let $Q\subseteq N(S)\setminus\{v\}$ be a small maximal nontrivial critical set.
We first show that no vertex of $Q\cap C$ has an owner.
By Lemma~\ref{res:route-3-24}, every vertex of $Q\cap C_+$ is joined to a vertex of $S\cap C_+$ by a thin edge.
Since $\alpha(S\cap C)=0$, both end vertices are incident with exactly one thin edge.
If a vertex of $Q\cap C_+$ had an owner, the other end vertex of its unique thin edge
would also have an owner by Lemma~\ref{route:lem:owner-pair}(iii), contrary to the fact that
no vertex of $S\cap C$ has an owner. Thus any vertex $z\in Q\cap C$ that has an owner must belong to
$C_0$ and be incident with exactly one thin edge. Its owner $p_0$ belongs to $Q$ by
Lemma~\ref{res:route-3-29}. Write $N(p_0)=\{a,b,r\}$ and $N(z)=\{a,b,w\}$, with
$za$ thin. As $za$ is the unique thin edge incident with $z$, the edge $zb$ is nonthin.
By Lemma~\ref{route:lem:owner-pair}(ii), $b\in P\cup C$,
$N(b)=\{p_0,z,z'\}$, and $rz'$ is the unique thin edge incident with $r\in C$.
Since $z\in C_0$, Lemma~\ref{res:route-3-28} gives $b\notin P$.
If $b\notin S$, then one of its neighbours $p_0,z$ must be $v$ by Lemma~\ref{res:route-3-13},
contrary to $p_0,z\in Q\subseteq N(S)\setminus\{v\}$. Hence $b\in S\cap C$. Since $bp_0,bz$ are
nonthin, $bz'$ is thin, and Lemma~\ref{res:route-3-28} gives $z'\in C$.
The two thin edges $bz',rz'$ then give
$\alpha(S\cap C)\ge\mu(z')\ge1$, a contradiction.
Therefore no vertex of $Q\cap C$ has an owner.

Put $m=|S\cap C_0|$ and $j=|N(S)\cap C_0|$.
Suppose now that $m+\frac j2+f_S<2$, so $m+f_S\le1$.
Let $R\subseteq S$ be a $2$-bounded set. Since no vertex of $P\cup D$ belongs to $R$,
write $R=X\cup V'$, where $X\subseteq C$ and $V'\subseteq V_{\ge4}$, and put $h=|V'|$.
Each vertex of $V'$ has degree at least four and is incident with at most two nonthin edges,
so it is incident with at least two thin edges.
By Lemma~\ref{res:route-3-8}(i), the other end vertices of these thin edges belong to $C$.
Since $e_S=0$, they belong to $N(S)\cap C_0$.
Each is incident with exactly one thin edge, so $j\ge2h$.
If $h\ge1$, the assumed inequality gives $h=1$ and $f_S=0$.
Since $|R|\ge2$, $X$ is nonempty.
Select two vertices joined by thin edges to the vertex of $V'$.
Neither belongs to $N(X)$, since an edge joining one of them to $X\subseteq C$
would be a nonthin edge counted by $f_S$, contrary to $f_S=0$.
Theorem~\ref{thm:brace-hall} now gives
\[
 |X|+3=|N(R)|\ge|N(X)|+2\ge|X|+4,
\]
a contradiction. Therefore $R\subseteq C$.

For each $z\in R\cap C_+$, let $zc_z$ be the unique thin edge incident with $z$.
Since $2\le|R|\le|S|\le b_q\le q-3$, the observation above gives another neighbour $z'\in R$ of $c_z$.
The vertices $c_z$ are distinct and are each incident with exactly one thin edge, so $c_zz'$ is nonthin.
These distinct edges are counted
by $f_S$, giving $f_S\ge|R\cap C_+|\ge2-m>f_S$, a contradiction.
Hence $S$ contains no $2$-bounded set.

By Lemma~\ref{comp:lem:root-deletion}, $|N(S\setminus\{u\})|=|S|+1$.
Also, $2\le|S|-1\le b_q-1\le q-3$. The set $S\setminus\{u\}$ contains no $2$-bounded set and,
since $S\in\mathcal A$, no vertex of $P$. Therefore it contains
$x_0\in U$ with at least three nonthin edges. By Lemma~\ref{res:route-3-26}, select two disjoint
small maximal nontrivial critical sets $Q_1,Q_2\subseteq N(S)\setminus\{v\}$
such that $|N(Q_i)\setminus S|\le2$ and $|Q_i|\le|S|-2$ for $i=1,2$.
No vertex of $(Q_1\cup Q_2)\cap C$ has an owner, by the first part of the proof.

For $i=1,2$, every vertex joined by a thin edge to a vertex of $Q_i\cap V_{\ge4}$ belongs to $S\cap C$
by Lemmas~\ref{res:route-3-8}(i) and~\ref{res:route-3-24}.
Its unique thin edge joins it to a vertex outside $C$, so it belongs to $S\cap C_0$.
Since $m\le1$, every vertex of $Q_i\cap V_{\ge4}$ is incident with at most one thin edge and hence
at least three nonthin edges.
As $S$ contains no $2$-bounded set, Lemma~\ref{res:route-3-33}, applied to $S$ and $Q_i$,
gives a $2$-bounded set $R_i\subseteq Q_i$.
No vertex of $P\cup D$ or $Q_i\cap V_{\ge4}$ can belong to $R_i$, since each is incident with at least
three nonthin edges. Thus $R_i\subseteq C$.
For $i=1,2$ and each $z\in R_i\cap C_+$, select $c_z\in C$ such that $zc_z$ is thin.
By Lemma~\ref{res:route-3-24}, $c_z\in S\cap C_+$.
Since $2\le|R_i|\le|Q_i|\le b_q\le q-3$, the observation above gives another neighbour $z'\in R_i$ of $c_z$.
The only thin edge incident with $c_z$ is $zc_z$, so $c_zz'$ is nonthin and is counted by $f_S$.
Different vertices $z$ give different vertices $c_z$, and hence different edges $c_zz'$.
The remaining vertices of $R_1\cup R_2$ belong to $N(S)\cap C_0$. Hence
$f_S+j\ge|R_1|+|R_2|\ge4$, and
$m+\frac j2+f_S\ge\frac{j+f_S}{2}\ge2$, the required contradiction.
\end{proof}

The following theorem implies Theorem~\ref{thm:main}.

\begin{theorem}\label{res:route-4-2}
Let $G$ be a simple brace of order $n\ge6$, and let $a(G)$ denote the number of cubic
vertices incident with three thin edges. Then $n+2\le6|T(G)|-2a(G)$.
\end{theorem}
\begin{proof}
Put $t=|T(G)|$. First suppose that $G$ is minimal and $n\ge12$.
Write $n=2q$, so $q\ge6$.
We first estimate the sums of the terms in Lemma~\ref{res:route-4-1} over all $S\in\mathcal A^*$.
Let $T_{CC}$ denote the set of thin edges with both end vertices in $C$, and let $T_{\ge4}$
denote the set of thin edges with one end vertex in $C$ and the other in $V_{\ge4}$.
By Lemma~\ref{res:route-3-8}(i), these sets partition $T(G)$.
Put $x=|T_{CC}|$ and $y=|T_{\ge4}|$, so $t=x+y$.

Let $g_+=|C_+|$ and $g_0=|C_0|$, and put $\xi=2x-g_+$ and $\eta=y-g_0$.
Also put $L=|E(P,C)|$ and $F=|E(G[C])\setminus T(G)|$, and let $\theta$
be the number of thin edges joining $C_+$ to $V_{\ge4}$. Put $\delta_0=\eta-\theta$.
Each edge of $T_{CC}$ is counted twice in $\sum_{c\in C_+}d_{T_{CC}}(c)$.
Of the $y$ edges in $T_{\ge4}$, exactly $\theta$ have an end vertex in $C_+$,
so the remaining $y-\theta$ have an end vertex in $C_0$. Hence
\begin{equation}
 \xi=\sum_{c\in C_+}(d_{T_{CC}}(c)-1),
\end{equation}
and
\begin{equation}\label{eq:route-4-1}
 \delta_0=\eta-\theta=\sum_{c\in C_0}(d_{T(G)}(c)-1).
\end{equation}
In particular, $\xi,\delta_0,\eta$ are nonnegative.
Also, by the definition of $\mu$,
$\sum_{c\in C}\mu(c)=2x+y-g_+-g_0=\xi+\eta$.

By Theorem~\ref{res:route-3-19}, each vertex outside $P$ belongs to at most one member of $\mathcal A$.
Let $c\in C_0$, and select a thin edge $cz$.
If $N(c)\cap S=\{z\}$ for some $S\in\mathcal A$, then $S$ is a nontrivial critical set
with root $z$ and designated neighbour $c$, contrary to the thinness of $cz$ by
Proposition~\ref{res:route-3-1} and Lemma~\ref{res:route-3-2}.
Thus every $S\in\mathcal A$ with $c\in N(S)$ contains at least one of the two neighbours of $c$ other than $z$.
These two neighbours lie outside $P$ by Lemma~\ref{res:route-3-28}, and each belongs to at most one
member of $\mathcal A$. Therefore $c\in N(S)$ for at most two members $S$ of $\mathcal A$.
Each vertex of $C_0$ also belongs to at most one member of $\mathcal A$.
By the definition of $w_S$ and since $\mu(c)\ge0$ for $c\in C_0$, it follows that
\[
 \sum_{S\in\mathcal A^*}w_S\le g_0+\frac12(2g_0)=2g_0,
\]
and
\[
 2\sum_{S\in\mathcal A^*}\sum_{c\in N(S)\cap C_0}\mu(c)
 \le4\sum_{c\in C_0}\mu(c)=4\delta_0.
\]
Every vertex of $V_{\ge4}$ belongs to at most one member of $\mathcal A$, so no edge is counted
by $e_S$ for two distinct sets $S\in\mathcal A^*$. All these edges join $C_+$ to $V_{\ge4}$.
Hence $\sum_{S\in\mathcal A^*}e_S\le\theta$.

We now consider the members of $\mathcal A$ with two vertices.
Put $N=|\mathcal A|$ and $a=a(G)$, and let $k$ be the number of components of $G[U]$.
By Lemma~\ref{res:route-3-21}(i), each such member is of the form $\{p_0,c\}$,
where $p_0\in P$ is an owner of $c\in C$.
Conversely, if $p_0$ is an owner of $c$, Lemma~\ref{res:route-3-29} gives
$\{p_0,c\}=S_{p_0,v}$, where $v$ is its designated neighbour, so $\{p_0,c\}\in\mathcal A$.
Let $N_1$ and $N_2$ denote the numbers of these sets for which $c$ is incident with exactly one
or exactly two thin edges, respectively.
By Lemma~\ref{route:lem:owner-pair}(i), each such vertex $c$ is incident with a nonthin edge
and hence with at most two thin edges.
By Lemma~\ref{res:route-3-29}, no vertex that has an owner belongs to a member of $\mathcal A^*$.

By Theorem~\ref{res:route-3-19}, each vertex $c\in C_+$ belongs to at most one member of $\mathcal A^*$.
Suppose that $c\in N_{T(G)}(S_i\cap C)$ for distinct $S_1,S_2\in\mathcal A^*$ and $i=1,2$.
For $i=1,2$, select $z_i\in S_i\cap C$ such that $cz_i$ is thin.
Each $S_i$ contains another neighbour of $c$, since otherwise it would be a nontrivial critical set
with root $z_i$ and designated neighbour $c$, contrary to the thinness of $cz_i$
by Proposition~\ref{res:route-3-1} and Lemma~\ref{res:route-3-2}.
Since $c$ is cubic, $S_1$ and $S_2$ contain a common neighbour of $c$.
By Theorem~\ref{res:route-3-19}, this vertex is their common root $p_0\in P$.
Lemma~\ref{res:route-3-28} then implies that $p_0$ is an owner of $z_1$.
This contradicts Lemma~\ref{res:route-3-29}, since $z_1\in S_1$ and $S_1\in\mathcal A^*$.

In $\sum_{S\in\mathcal A^*}\alpha(S\cap C)$, the term $\mu(c)$ is therefore counted at most twice for $c\in C_+$.
For $c\in C_0$, it is counted only when $c\in S$, and hence at most once.
By Lemma~\ref{res:route-3-29}, there are $N_2$ distinct vertices $c$ that have an owner
and are incident with exactly two thin edges.
For each, $\mu(c)=1$ and $c$ belongs to no member of $\mathcal A^*$.
Its term is counted at most once if $c\in C_+$, and not at all if $c\in C_0$, decreasing the bound by one.
Each of the $a$ vertices incident with three thin edges has $\mu(c)=2$ and $\varepsilon(c)=1$.
If it belongs to a member of $\mathcal A^*$, subtracting $\varepsilon(c)$ decreases the bound by one;
otherwise, its term $\mu(c)=2$ for membership in $S$ is absent, decreasing the bound by at least two.
By Lemma~\ref{route:lem:owner-pair}(i), these $a$ vertices are different from the $N_2$ vertices.
Since $\sum_{c\in C_+}\mu(c)=\xi+\theta$ and $\sum_{c\in C_0}\mu(c)=\delta_0$, we obtain
\[
 \sum_{S\in\mathcal A^*}\left(\alpha(S\cap C)
       -\sum_{c\in S\cap C}\varepsilon(c)\right)
 \le2(\xi+\theta)+\delta_0-N_2-a
 =2\xi+\eta+\theta-N_2-a.
\]

Let $b_U,b_P,b_C$ denote the numbers of large sets $S_{u,v}$ with $v\in U$, $v\in P$ and $v\in C$,
respectively, so $\ell=b_U+b_P+b_C$.
For each of the $b_P$ sets $S_{p_0,b}$, take $c_1,c_2$ and $w_1,w_2$ as in Lemma~\ref{res:route-3-37}.
Let $\ell_C$ denote the total number of these vertices $c_i$ for which $w_i\in C$.
These pairs $(c_i,c_iw_i)$ are distinct and $c_iw_i$ is nonthin by Lemma~\ref{res:route-3-37}(ii).
Its proof also shows that each $c_i$ belongs to $C$ and $p_0$ is an owner of $c_i$.
For each edge $e$ counted by $f_S$, consider the pair $(c,e)$, where $c$ is its unique end vertex in $S\cap C$.
By Theorem~\ref{res:route-3-19}, the same pair cannot be counted for two distinct members of $\mathcal A^*$.
There are $2F$ pairs $(c,e)$ with $e$ a nonthin edge of $G[C]$ and $c$ an end vertex of $e$.
Of these, the $\ell_C$ pairs $(c_i,c_iw_i)$ cannot be counted by any $f_S$,
since $c_i$ belongs to no member of $\mathcal A^*$ by Lemma~\ref{res:route-3-29}.
Hence $\sum_{S\in\mathcal A^*}f_S\le2F-\ell_C$.
Summing the inequality in Lemma~\ref{res:route-4-1} over all $S\in\mathcal A^*$,
applying these bounds and using $\delta_0=\eta-\theta$, we obtain
\[
 2|\mathcal A^*|\le2g_0+4\xi+6\eta+2F-2a-2N_2-\ell_C.
\]

For each set $\{p_0,c\}$ counted by $N_1$, let $v$ be its designated neighbour.
By Lemma~\ref{route:lem:owner-pair}(ii), $v\in C$, so $p_0v\in E(P,C)$.
Lemma~\ref{res:route-3-29} implies that distinct sets counted by $N_1$ give distinct edges $p_0v$.
The $b_C$ large sets $S_{u,v}$ also give distinct edges $uv\in E(P,C)$,
since only one set was selected for each root and designated neighbour.
None is an edge obtained from a set counted by $N_1$, since Lemma~\ref{res:route-3-29}
makes $\{p_0,c\}$ the only nontrivial critical set with root $p_0$ and designated neighbour $v$.
Thus $N_1+b_C\le L$.
Put $\rho=2g_0+4\xi+6\eta+2L+2F$.
Since $N=|\mathcal A^*|+N_1+N_2$, the preceding inequalities give $2N+2a+2b_C+\ell_C\le\rho$.

Let $\ell_U$ denote the number of vertices $c_i$ obtained above for which $w_i\in U$.
By Lemma~\ref{res:route-3-37}, each of the $b_P$ sets gives two distinct vertices $c_i$,
all these vertices are distinct, and every $w_i$ belongs to $U\cup C$. Hence $\ell_U+\ell_C=2b_P$.
By Lemma~\ref{res:route-3-11}(iii), at most one of the sets $S_{u,v}$ with a given root $u\in P$ is large.
Thus the $b_U$ sets have distinct roots in $P$, all different from the $\ell_U$ distinct roots $c_i\in C$.
For each of these $b_U+\ell_U$ roots, there is a large nontrivial critical set whose designated neighbour belongs to $U$.
By Lemma~\ref{res:route-3-36}, these designated neighbours belong to distinct components of $G[U]$.
Hence $b_U+\ell_U\le k$, and consequently $2b_U+2b_P\le2k+\ell_C$.
Together with $\ell=b_U+b_P+b_C$ and $2N+2a+2b_C+\ell_C\le\rho$, this gives
\begin{equation}\label{eq:route-4-2}
 2N+2\ell+2a\le\rho+2k.
\end{equation}

Put $E=|E(G[P])|$. By Lemma~\ref{res:route-3-8}(iii), $G[U]$ is a forest with $k$ components
and hence has $|U|-k$ edges. Since every vertex outside $U$ is cubic,
$\sum_{u\in U}d(u)=3|U|+\sigma$. Thus
$|E(U,C\cup P)|=3|U|+\sigma-2(|U|-k)=|U|+\sigma+2k$.
On the other hand, every vertex of $C\cup P$ is cubic, and $G[C\cup P]$ has $x+F+E+L$ edges.
Counting the edges joining $U$ to $C\cup P$ at their end vertices in $C\cup P$, we obtain
\[
 |U|+\sigma+2k=3(g_++g_0+p)-2(x+F)-2E-2L.
\]
Substituting $n=|U|+g_++g_0+p$, $g_+=2x-\xi$ and $g_0=y-\eta$ and using
$t=x+y$, we obtain
\begin{equation}\label{eq:route-4-3}
 6t-n=\rho+2E+\sigma+2k-4p.
\end{equation}

If $p>0$, \eqref{eq:route-4-3}, \eqref{eq:route-4-2} and Theorem~\ref{res:route-3-19} give
\[
 6t-n-2a\ge2E+\sigma+2N+2\ell-4p\ge2E+2\ge2.
\]
If $p=0$, then $P=\varnothing$ and $E=L=0$.
Also, $\xi+\eta=\sum_{c\in C}\mu(c)\ge2a$, so \eqref{eq:route-4-3} gives
\[
\begin{aligned}
 6t-n-2a
 &=2g_0+4\xi+6\eta+2F+\sigma+2k-2a\\
 &\ge2g_0+3\xi+5\eta+2F+\sigma+2k\ge0.
\end{aligned}
\]
Suppose that $6t-n-2a=0$. Since every term in the last sum is nonnegative,
$g_0=\xi=\eta=F=\sigma=k=0$, and $2a\le\xi+\eta$ gives $a=0$.
As $k=0$, we have $U=\varnothing$. Together with $p=g_0=0$, this implies that every vertex lies in $C_+$.
Since $\xi=0$, each vertex is incident with exactly one thin edge.
Its other two incident edges are nonthin and have both end vertices in $C$, contrary to $F=0$.
Thus $6t-n-2a>0$. Since $n=2q$, the integer $6t-n-2a$ is even, and hence $6t-n-2a\ge2$.

Now let $G$ be an arbitrary simple brace of order $n\ge12$.
By Lemma~\ref{prop:spanning-minimal-brace}, select a spanning minimal brace $H$ of $G$.
Proposition~\ref{res:route-3-4} gives $T(H)\subseteq T(G)$.
Let $v$ be a cubic vertex of $G$ incident with three thin edges.
Since $H$ is spanning and has minimum degree at least three, $v$ has the same three incident edges in $H$ as in $G$.
If $v$ is not incident with three thin edges in $H$, then at least one of its incident edges belongs to $T(G)\setminus T(H)$.
Each edge of $T(G)\setminus T(H)$ has two end vertices, so at most
$2\bigl(|T(G)|-|T(H)|\bigr)$ of the vertices counted by $a(G)$ are not counted by $a(H)$.
Hence $a(G)-a(H)\le2\bigl(|T(G)|-|T(H)|\bigr)$.
Applying the inequality already proved for $H$, we obtain
\[
 6|T(G)|-2a(G)\ge6|T(H)|-2a(H)+2\bigl(|T(G)|-|T(H)|\bigr)\ge n+2.
\]

Finally, suppose that $n=6,8$ or $10$. By~\cite{CLM}, $G$ has at least two thin edges, so $t\ge2$.
If $t=2$, no vertex is incident with three thin edges, so $a(G)=0$ and $6t-2a(G)=12\ge n+2$.
If $t\ge3$, counting thin edges at their end vertices gives $3a(G)\le2t$.
Thus $6t-2a(G)\ge\frac{14}{3}t\ge14>n+2$. This proves the theorem.
\end{proof}

\begin{proof}[\normalfont\bfseries Proof of Theorem~\ref{thm:main}]
Let $G$ be a simple brace of order $n\ge6$. Since $a(G)\ge0$,
Theorem~\ref{res:route-4-2} gives $6|T(G)|\ge n+2$. As $|T(G)|$ is an integer,
$|T(G)|\ge\lceil\frac{n+2}{6}\rceil$.
\end{proof}

\subsection{The Proof of Theorem~\ref{thm:general}}\label{sec:proof-general}
We begin with the following lemma.

\begin{lemma}\label{lem:underlying-simple}
Let $G$ be a graph, and let $G_1$ be its underlying simple graph.
Then $G$ is matching covered if and only if $G_1$ is matching covered.
Moreover, $G$ is a brace if and only if $G_1$ is a brace.
\end{lemma}
\begin{proof}
Replacing each edge of a perfect matching of $G$ by the edge of $G_1$ with the same end vertices
gives a perfect matching of $G_1$. Conversely, for each edge of a perfect matching of $G_1$,
select an edge of $G$ with the same end vertices. The selected edges form a perfect matching of $G$.
Thus an edge of $G$ belongs to a perfect matching if and only if the edge of $G_1$ with the same
end vertices belongs to a perfect matching.
Since $G$ is connected if and only if $G_1$ is connected, the first assertion follows.
Now suppose that $G$ and $G_1$ are matching covered. The corresponding perfect matchings meet
$\partial_G(X)$ and $\partial_{G_1}(X)$ in the same number of edges for every
nonempty proper subset $X$ of $V(G)$. Thus $\partial_G(X)$ is tight if and only if
$\partial_{G_1}(X)$ is tight. Since $V(G)=V(G_1)$, the two cuts are either both nontrivial
or both trivial. Also, $G$ is bipartite if and only if $G_1$ is bipartite.
The second assertion follows.
\end{proof}

Here, we give the proof of Theorem~\ref{thm:general}.

\begin{proof}[\normalfont\bfseries Proof of Theorem~\ref{thm:general}]
Let $G$ be a brace of order $n\ge6$,
and let $G_1$ be its underlying simple graph. By Lemma~\ref{lem:underlying-simple}, $G_1$ is a brace.
Let $\mathcal D$ be the set of edges of $G_1$ whose end vertices are joined by at least two edges of $G$, and
let $C_{\mathcal D}$ be the set of cubic vertices of $G_1$ incident with an edge of $\mathcal D$.
For $uv\in E(G_1)\setminus\mathcal D$, we also use $uv$ to denote the unique edge of $G$ joining $u$ and $v$.
Let $e$ be an edge of $G$ that has a parallel edge. The underlying simple graph of $G-e$ is $G_1$.
By Lemma~\ref{lem:underlying-simple}, $G-e$ is a brace. Since its order is at least six,
Theorem~\ref{thm:brace-connectivity} implies that it has minimum degree at least three.
Thus $G-e$ is its own retract, and $e$ is thin.

Next, let $e\in T(G_1)\setminus\mathcal D$ have no end vertex in $C_{\mathcal D}$.
Adding parallel edges to $G_1-e$ gives $G-e$, which is matching covered by Lemma~\ref{lem:underlying-simple}.
By Lemma~\ref{lem:edge-retract}, the retract of $G_1-e$ is obtained by
bicontracting the cubic end vertices of $e$ in $G_1$. Each such end vertex $x$ lies outside
$C_{\mathcal D}$, so all edges of $G$ incident with $x$ have multiplicity one.
Thus $x$ has degree two and the same two neighbours in $G_1-e$ and $G-e$.
By the same lemma, the sets $\{x\}\cup N_{G_1-e}(x)$ contracted at these vertices are disjoint.
Since $G$ is bipartite, the only edges within each such set are the two edges incident with $x$.
The same bicontractions in $G-e$ therefore give the retract of $G_1-e$ with parallel edges added.
By the proof of Lemma~\ref{lem:edge-retract}, this retract has minimum degree at least three.
Adding parallel edges cannot decrease any vertex degree, so no further bicontraction is required.
Since $e$ is thin in $G_1$, the retract of $G_1-e$ is a brace. By Lemma~\ref{lem:underlying-simple},
the resulting graph is also a brace, so $e$ is thin in $G$.

Let $d=|\mathcal D|$.
For each $c\in C_{\mathcal D}$, select an edge $f(c)\in\mathcal D$ incident with $c$.
Let $s$ count the edges of $T(G_1)\setminus\mathcal D$ having no end vertex in $C_{\mathcal D}$.
These $s$ edges are thin in $G$.
For each $uv\in\mathcal D$, at least two edges of $G$ join $u$ and $v$, and all of them are thin.
Hence
\begin{equation}\label{eq:route-5-1}
 |T(G)|\ge s+2d.
\end{equation}
For $c\in C_{\mathcal D}$, put $r_c=d_{T(G_1)}(c)$ and $A_c=6r_c-2\mathbf1_{\{r_c=3\}}$, where the
indicator is one when $r_c=3$ and zero otherwise. Since $f(c)\in\mathcal D$, at most
$r_c-\mathbf1_{\{f(c)\in T(G_1)\}}$ edges of $T(G_1)\setminus\mathcal D$ are incident with $c$.
Every edge of $T(G_1)\setminus\mathcal D$ not counted by $s$ has an end vertex in $C_{\mathcal D}$.
Hence
\[
 |T(G_1)|\le s+|\mathcal D\cap T(G_1)|+\sum_{c\in C_{\mathcal D}}\bigl(r_c-\mathbf1_{\{f(c)\in T(G_1)\}}\bigr).
\]
By Theorem~\ref{res:route-4-2}, $n\le6|T(G_1)|-2a(G_1)$. Also $a(G_1)\ge\sum_{c\in C_{\mathcal D}}\mathbf1_{\{r_c=3\}}$,
since every vertex of $C_{\mathcal D}$ is cubic. Substitution gives
\[
 n\le6s+6|\mathcal D\cap T(G_1)|+\sum_{c\in C_{\mathcal D}}\bigl(A_c-6\mathbf1_{\{f(c)\in T(G_1)\}}\bigr).
\]
Group the summands according to the edge $f(c)$. For each $f\in\mathcal D$, at most two
vertices $c\in C_{\mathcal D}$ satisfy $f(c)=f$, since each such vertex is an end vertex of $f$.
If $f(c)$ is thin in $G_1$,
then $1\le r_c\le3$ and $A_c\le16$. Thus an edge of $\mathcal D\cap T(G_1)$, including its
contribution to $6|\mathcal D\cap T(G_1)|$, accounts for at most $6+2(16-6)=26$.
If $f(c)$ is nonthin in $G_1$, then $r_c\le2$ and $A_c\le12$, so an edge of $\mathcal D\setminus T(G_1)$
accounts for at most $24$. Therefore
\[
 n\le6s+26|\mathcal D\cap T(G_1)|+24|\mathcal D\setminus T(G_1)|\le13(s+2d)\le13|T(G)|,
\]
where the last inequality is \eqref{eq:route-5-1}.
Since $|T(G)|$ is an integer, $|T(G)|\ge\lceil\frac{n}{13}\rceil$.
\end{proof}

\section{Thin Edges in Cubic Braces}\label{sec:cubic}
In this section, we prove Theorem~\ref{cub:thm:main}.

\subsection{The Number of Thin Edges}
We begin with some notation.
Let $G$ be a cubic brace of order $2q$, where $q\ge6$.
The sets $C,D,P$ are those defined in Subsection~3.2.
Let $C_{PP}$ consist of the vertices of $C$ adjacent to an end vertex of an edge of $G[P]$.
Let $I$ consist of the thin edges whose end vertices lie outside $C_{PP}$ and are each incident with exactly one thin edge.
Let $V(I)$ be the set of their end vertices.

\begin{lemma}\label{lem:cubic-owners}
Let $G$ be a cubic brace of order at least twelve. The following statements hold.
\begin{enumerate}[label=\textup{(\roman*)},leftmargin=*]
\item A vertex $c\in C$ has an owner and satisfies $d_{T(G)}(c)=1$ if and only if $c\in C_{PP}$.
\item If $z$ has an owner and $d_{T(G)}(z)=2$, then every vertex joined to $z$ by a thin edge
lies outside $C_{PP}$ and is incident with exactly one thin edge.
\item Every vertex of $C$ adjacent to a vertex of $P$ is incident with exactly one thin edge.
No vertex of $V(I)$ is adjacent to a vertex of $P$.
\end{enumerate}
\end{lemma}
\begin{proof}
By Theorem~\ref{thm:brace-connectivity}, each vertex of $G$ has three distinct neighbours, so $G$ is simple.
Deleting an edge leaves vertices of degree two, so the same theorem implies that $G$ is minimal.
Since $G$ is cubic, $C=C_+$.
By Lemma~\ref{res:route-3-34}, the neighbours in $C$ of the two end vertices of each edge of $G[P]$
are precisely the four corners of the corresponding domino subgraph.
Thus $C_{PP}$ is the set of all these corners.
By Lemma~\ref{route:lem:owner-pair}(i), a vertex having an owner is incident with one or two thin edges.
Every such vertex incident with one thin edge belongs to $C_{PP}$ by Lemma~\ref{res:route-3-34}(ii).
Conversely, each corner in a subgraph given by Lemma~\ref{res:route-3-34} is incident with exactly one thin edge.
The other end vertex of that edge is adjacent to a centre in $P$,
so Lemma~\ref{res:route-3-28} gives an owner of the corner. This proves (i).

For (ii), let $u$ be an owner of $z$, and write $N(u)=\{a,b,v\}$ and $N(z)=\{a,b,w\}$.
By Lemma~\ref{route:lem:owner-pair}(i), $zw$ is nonthin, so $za,zb$ are thin.
Write $N(a)=\{u,z,a'\}$ and $N(b)=\{u,z,b'\}$.
Theorem~\ref{thm:brace-hall} applied to $\{a,b\}$ gives $a'\ne b'$.
Thus $\{a,b\}$ is a nontrivial critical set with root $a$ and designated neighbour $a'$,
and also with root $b$ and designated neighbour $b'$.
By Proposition~\ref{res:route-3-1} and Lemma~\ref{res:route-3-2}, $aa',bb'$ are nonthin.
Since $au,bu$ are nonthin, $az,bz$ are the unique thin edges incident with $a,b$, respectively.
If either $a$ or $b$ belonged to $C_{PP}$, its unique thin edge would join it to another vertex of $C_{PP}$.
That vertex would be $z$, contrary to $d_{T(G)}(z)=2$. This proves (ii).

For (iii), let $c\in C$ be adjacent to $u\in P$, and select a thin edge $cz$.
By Lemma~\ref{res:route-3-28}, $u$ is an owner of $z$.
If $d_{T(G)}(z)=1$, then (i) gives $z\in C_{PP}$, and its unique thin edge also has its other end vertex in $C_{PP}$.
If $d_{T(G)}(z)=2$, apply (ii).
In either case, $c$ is incident with exactly one thin edge.
If an end vertex $c$ of an edge $cz\in I$ had a neighbour in $P$,
Lemma~\ref{res:route-3-28} would give an owner of $z$.
Since $d_{T(G)}(z)=1$, part (i) would imply $z\in C_{PP}$, contrary to the definition of $I$.
\end{proof}

We next consider the family $\mathcal A$ from Subsection~3.3.
Its members are the small maximal nontrivial critical sets $S_{u,v}$ with $S_{u,v}\cap P=\{u\}$.
Put $\mathcal F=E(G[C])\setminus T(G)$.
For $c\in V(I)$, let $c^*$ denote the other end vertex of the unique thin edge incident with $c$.
For each edge $cc^*\in I$, put $w(cc^*)=\max\{0,d_{\mathcal F}(c)+d_{\mathcal F}(c^*)-1\}$.
For $X\subseteq C$ contained in one part of the bipartition, write
$\alpha(X)=\sum_{z\in C\cap(X\cup N_{T(G)}(X))}(d_{T(G)}(z)-1)$,
as in Subsection~\ref{sec:proof-main}, and put
$\omega(X)=\sum_{c\in X\cap V(I)}w(cc^*)$.

\begin{lemma}\label{lem:cubic-small-set}
Let $G$ be a cubic brace of order $2q$, where $q\ge6$.
For every $S\in\mathcal A$ with $|S|\ge3$,
\[
 4\le\omega(S\cap C)+2\alpha(S\cap C)+
 2\sum_{c\in(S\cap C)\setminus V(I)}d_{\mathcal F}(c).
\]
\end{lemma}
\begin{proof}
Since $G$ is cubic, $U=D$. We use $b_q=\lfloor\frac{q-2}{2}\rfloor$, as in Subsection~3.2.
Let $X\subseteq V(I)$ be contained in a part of the bipartition, with $|N(X)|=|X|+2$ and $2\le|X|\le q-3$.
We first show that $\omega(X)\ge|X|$.
For each $c\in X$, the vertex $c^*$ has
a neighbour in $X$ other than $c$; otherwise $X$ would be a nontrivial critical
set with root $c$ and designated neighbour $c^*$.
Proposition~\ref{res:route-3-1} and Lemma~\ref{res:route-3-2} would then imply that $cc^*$ is nonthin, contrary to $cc^*\in I$.
For each $c\in X$, select one edge joining $c^*$ to $X\setminus\{c\}$.
This edge belongs to $\mathcal F$, since both end vertices lie in $C$ and $cc^*$ is the unique thin edge incident with $c^*$.
Since $I$ is a matching, the vertices $c^*$ are distinct and lie in the part not containing $X$.
The selected edges are therefore distinct. Counting these $|X|$ edges at the vertices $c^*$ and at their end vertices in $X$ gives
\[
 \sum_{c\in X}d_{\mathcal F}(c^*)\ge|X|,
\]
and
\[
 \sum_{c\in X}d_{\mathcal F}(c)\ge|X|.
\]
Thus $\omega(X)\ge\sum_{c\in X}(d_{\mathcal F}(c)+d_{\mathcal F}(c^*)-1)\ge|X|$.

Now let $S\in\mathcal A$ have at least three vertices, with root $u$ and designated neighbour $v$.
By Lemma~\ref{res:route-3-29}, no vertex of $S\cap C$ has an owner.
In particular, $S\cap C_{PP}=\varnothing$ by Lemma~\ref{lem:cubic-owners}(i).
If $\alpha(S\cap C)\ge2$, the inequality follows. Assume henceforth that $\alpha(S\cap C)\le1$.

Suppose first that $S\setminus\{u\}$ contains a vertex $x\in D$.
Then $x\in U$, and all three edges incident with $x$ are nonthin.
By Lemma~\ref{res:route-3-26}(i), (ii), select two disjoint small maximal nontrivial critical sets
$Q_1,Q_2\subseteq N(S)\setminus\{v\}$ such that $|N(Q_i)\setminus S|\le2$
and $|Q_i|\le|S|-2$ for $i=1,2$. By Lemma~\ref{res:route-3-24}, every vertex joined by a thin edge to a vertex of $Q_1\cup Q_2$ belongs to $S$.
Thus each vertex of $Q_i\cap C$ is joined by a thin edge to a vertex of $S\cap C$,
and all vertices joined to $Q_i\cap C$ by thin edges belong to $S\cap C$.
Consequently $\alpha(Q_i\cap C)\le\alpha(S\cap C)$ for $i=1,2$.
If $\alpha(S\cap C)=1$, its defining sum has a unique nonzero term,
corresponding to a vertex $z$ with $d_{T(G)}(z)=2$.
If $z\in Q_1\cup Q_2$, it belongs to only one of these disjoint sets and has no neighbour in either set,
since $Q_1,Q_2$ lie in the same part of the bipartition.
Otherwise, for $z$ to occur in both sums, it must be joined by a thin edge to a vertex of each $Q_i$.
Each $Q_i$ would then contain a second neighbour of $z$, since otherwise that thin edge would be nonthin
by Proposition~\ref{res:route-3-1} and Lemma~\ref{res:route-3-2}.
This would give four distinct neighbours of the cubic vertex $z$.
Hence at most one of $\alpha(Q_1\cap C),\alpha(Q_2\cap C)$ is positive.

Select $i\in\{1,2\}$ with $\alpha(Q_i\cap C)=0$.
Both end vertices of every thin edge incident with $Q_i\cap C$ are incident with exactly one thin edge.
If one end vertex belonged to $C_{PP}$, both would belong to $C_{PP}$ by Lemma~\ref{res:route-3-34}(i),
contrary to $S\cap C_{PP}=\varnothing$, since one end vertex lies in $S$.
The edge therefore belongs to $I$, and neither end vertex has an owner by Lemma~\ref{lem:cubic-owners}(i).
If $Q_i$ contains a $2$-bounded set $X$, every vertex of $X$ is cubic and incident with at most two nonthin edges.
Hence $X\subseteq Q_i\cap C\subseteq V(I)$.
Also, $|N(X)|=|X|+2$ and $2\le|X|\le|Q_i|\le q-3$, so the first paragraph gives $\omega(X)\ge|X|\ge2$.
For each $c\in X$, the vertex $c^*$ belongs to $S\cap C$.
These vertices $c^*$ are distinct, since $I$ is a matching, and $w(cc^*)=w(c^*c)$.
Thus all the terms of $\omega(X)$ occur in $\omega(S\cap C)$ without repetition.

If $Q_i$ contains no $2$-bounded set, Lemma~\ref{res:route-3-32} gives two disjoint small maximal
nontrivial critical sets contained in $N(Q_i)$. Since they are disjoint, at least one does not contain
the designated neighbour of $Q_i$. Select such a set $R$.
Then $|N(R)\setminus Q_i|\le2$ and $|R|\le|Q_i|-2$.
By Lemma~\ref{res:route-3-24}, every vertex joined by a thin edge to a vertex of $R\cap C$ belongs to $Q_i$.
Each such edge is incident with $Q_i\cap C$ and therefore belongs to $I$; neither end vertex has an owner.
Applying Lemma~\ref{res:route-3-33} to $Q_i$ and $R$ gives a $2$-bounded set $X\subseteq R$,
since $Q_i$ contains no such set. As $G$ is cubic, every vertex of $X$ is incident with a thin edge,
so $X\subseteq R\cap C\subseteq V(I)$. Also, $|N(X)|=|X|+2$ and $2\le|X|\le|R|\le q-3$.
The first paragraph therefore gives $\omega(X)\ge|X|\ge2$.
For each $c\in X$, the vertex $c^*$ belongs to $Q_i$, so Lemma~\ref{res:route-3-24}, applied to $S$ and $Q_i$,
gives $c\in S$. Thus $X\subseteq S$.

In either case, the summands $w(cc^*)$ in $\omega(S\cap C)$ for edges
$cc^*$ meeting $Q_i$ have sum at least two. The edges meeting $Q_1$ and those
meeting $Q_2$ are distinct, since these sets are disjoint and lie in the same part.
If $\alpha(Q_i\cap C)=0$ for both $i$, then $\omega(S\cap C)\ge4$.
Otherwise one of these two values is positive, so $1\le\alpha(S\cap C)\le1$.
Together with $\omega(S\cap C)\ge2$, this proves the inequality when $S\setminus\{u\}$ contains a vertex of $D$.

Suppose next that $S\setminus\{u\}$ contains no vertex of $D$.
Since $S\cap P=\{u\}$, we may write $S=\{u\}\cup X$ with $X\subseteq C$.
Put $m=|X|$, so $2\le m\le b_q-1\le q-3$. Lemma~\ref{comp:lem:root-deletion} gives $|N(X)|=|X|+2$.
Write $N(u)=\{v,a,b\}$ and select nontrivial critical sets $S_a,S_b$ with root $u$
and respective designated neighbours $a,b$.
By Lemma~\ref{comp:lem:root-deletion}, $a$ has a neighbour in each of $X$ and $S_b\setminus\{u\}$.
These sets are disjoint by Lemma~\ref{res:route-3-10}(i). Since $a$ is also adjacent to $u$ and $d(a)=3$,
it has exactly one neighbour in $X$. Applying the same argument to $b$ and $S_a$ gives exactly one neighbour of $b$ in $X$.
If a thin edge joined $a$ or $b$ to a vertex of $X$, Lemma~\ref{res:route-3-28} would give $u$ as an owner of that vertex,
contrary to the fact that no vertex of $S\cap C$ has an owner. Hence neither $a$ nor $b$ belongs to $N_{T(G)}(X)$.

If $m=2$, write $X=\{c,d\}$. If $a,b$ had the same neighbour in $X$, say $c$,
then $N(d)$ would be contained in the two-element set $N(S)\setminus\{v,a,b\}$, contrary to $d(d)=3$.
Thus $a,b$ have distinct neighbours in $X$.
Since $|N(S)\setminus\{v\}|=4$, we may write
$N(c)=\{a,x,y\}$ and $N(d)=\{b,x,y\}$. If any of $cx,cy,dx,dy$ is nonthin, adjust notation so that it is $cx$.
A nontrivial critical set contained in $S$, with root $c$ and designated neighbour $x$, cannot contain $d$,
since $xd\in E(G)$. Such a set would therefore be $\{u,c\}$, but $|N(\{u,c\})|=5$, not four.
Thus $S$ contains no such critical set. Lemma~\ref{res:route-3-23} therefore
gives a small maximal nontrivial critical set $Q$ with root $x$ and designated neighbour $c$ such that
$Q\subseteq N(S)$ and $|N(Q)\setminus S|\le2$.
Since $N(c)\cap Q=\{x\}$ and $N(S)=\{v,a,b,x,y\}$, neither $a$ nor $y$ belongs to $Q$.
Hence $x\in Q\subseteq\{x,b,v\}$.

If $v\in Q$, write $N(v)=\{u,r,s\}$. Since $r,s\notin S$,
we have $N(Q)\setminus S=\{r,s\}$. By Lemma~\ref{comp:lem:root-deletion},
$v$ has a neighbour in each of $S_a\setminus\{u\}$ and $S_b\setminus\{u\}$;
these two sets are disjoint by Lemma~\ref{res:route-3-10}(i).
After interchanging $r,s$ if necessary, we have $r\in S_a\setminus\{u\}$ and $s\in S_b\setminus\{u\}$. By Lemma~\ref{res:route-3-11}(ii),
$x\in N(S)\setminus N(u)$ cannot be adjacent to either $r$ or $s$.
Thus $N(x)\subseteq S\setminus\{u\}=\{c,d\}$, contrary to $d(x)=3$.
Hence $v\notin Q$, so $Q=\{x,b\}$. Lemma~\ref{res:route-3-21}(i) gives $b\in C$.
Select a thin edge $bz$. By Lemma~\ref{res:route-3-28}, $u$ is an owner of $z$.
Since $b\in Q\subseteq N(S)\setminus\{v\}$, Lemma~\ref{res:route-3-24} gives $z\in S$,
contrary to the fact that no vertex of $S\cap C$ has an owner.
Hence $cx,cy,dx,dy$ are thin, and each of $c,d$ is incident with at least two thin edges.
Thus $\alpha(X)\ge2$, contrary to $\alpha(S\cap C)\le1$, since $X=S\cap C$.

We may now suppose $m\ge3$ and $\alpha(X)\le1$.
Since $X\subseteq S$ and $S\cap C_{PP}=\varnothing$, Lemma~\ref{res:route-3-34}(i) implies that
no thin edge incident with $X$ has an end vertex in $C_{PP}$.
Let $c\in X\cap V(I)$. If $N(c^*)\cap X=\{c\}$, then $X$ would be a nontrivial critical set
with root $c$ and designated neighbour $c^*$, contrary to the thinness of $cc^*$.
Hence $c^*$ has a neighbour in $X\setminus\{c\}$. The edge joining them belongs to $\mathcal F$,
since $cc^*$ is the unique thin edge incident with $c^*$.
Thus $d_{\mathcal F}(c^*)\ge1$ and $w(cc^*)\ge d_{\mathcal F}(c)$.
If $\alpha(X)=0$, both end vertices of every thin edge incident with $X$ are incident with exactly one thin edge,
so $X\subseteq V(I)$. The $m$ vertices $c^*$, for $c\in X$, are distinct and different from $a,b$.
Together with $a,b$, they form $N(X)$, since $|N(X)|=m+2$.
Of the $3m$ edges incident with vertices of $X$, exactly $m$ are thin and two join $X$ to $a,b$.
The remaining $2m-2$ edges belong to $\mathcal F$.
Therefore $\omega(X)\ge\sum_{c\in X}d_{\mathcal F}(c)\ge2m-2\ge4$.

If $\alpha(X)=1$, exactly one vertex $z\in X\cup N_{T(G)}(X)$ is incident with more than one thin edge,
and $d_{T(G)}(z)=2$. If $z\in X$, exactly $m+1$ thin edges are incident with vertices of $X$.
Their end vertices outside $X$ are distinct, since each is incident with exactly one thin edge.
Together with $a,b$, they give $m+3$ distinct vertices of $N(X)$, contrary to $|N(X)|=m+2$.
Hence $z\notin X$. At most two vertices of $X$ are joined to $z$ by thin edges.
For every other vertex of $X$, both end vertices of its unique thin edge are incident with exactly one thin edge,
so that edge belongs to $I$. Since $m\ge3$, $X\cap V(I)$ is nonempty.
If an edge of $\mathcal F$ is incident with $X\setminus V(I)$, then
$\sum_{c\in X\setminus V(I)}d_{\mathcal F}(c)\ge1$, which, together with $2\alpha(X)=2$, proves the inequality.
Suppose that no edge of $\mathcal F$ is incident with $X\setminus V(I)$.
For each $c\in X\cap V(I)$, select a neighbour $c'\in X\setminus\{c\}$ of $c^*$.
The edge $c^*c'$ belongs to $\mathcal F$, so $c'\in X\cap V(I)$. In particular, $|X\cap V(I)|\ge2$.
Since $I$ is a matching, the vertices $c^*$ are distinct and lie in the same part of the bipartition.
The selected edges are therefore distinct. Counting them at their end vertices $c^*$ and $c'$ gives
$\sum_{c\in X\cap V(I)}d_{\mathcal F}(c^*)\ge|X\cap V(I)|$ and
$\sum_{c\in X\cap V(I)}d_{\mathcal F}(c)\ge|X\cap V(I)|$. Consequently
$\omega(X)\ge|X\cap V(I)|\ge2$, and $\omega(X)+2\alpha(X)\ge4$.
This completes the proof.

\end{proof}

We now sum this estimate to obtain the bound for cubic braces.

\begin{lemma}\label{res:route-6-2}
Every cubic brace $G$ of order at least twelve satisfies $5|T(G)|>|V(G)|$.
\end{lemma}
\begin{proof}
Let $G$ have order $n=2q$, and put $t=|T(G)|$, $p=|P|$, $E=|E(G[P])|$,
$L=|E(P,C)|$ and $N=|\mathcal A|$.
Let $k$ be the number of components of $G[D]$,
and let $M$ be the number of members of $\mathcal A$ with at least three vertices.
Let $N_2$ be the number of vertices that have an owner and are incident with exactly two thin edges.
By Lemma~\ref{res:route-3-34}, the $E$ domino subgraphs are vertex-disjoint and each has four corners, so $|C_{PP}|=4E$.
By Lemma~\ref{lem:cubic-owners}(i), exactly $4E$ vertices have an owner and are incident with one thin edge.
Each contributes a distinct edge to $E(P,C)$ in the subgraphs given by Lemma~\ref{res:route-3-34}.
If $z$ is counted by $N_2$ and has owner $u$, the two vertices joined to $z$ by thin edges are neighbours of $u$.
By Lemma~\ref{lem:cubic-owners}(ii), these vertices lie outside $C_{PP}$ and are each incident with one thin edge.
They are distinct for different vertices $z$ counted by $N_2$.
Hence $L\ge4E+2N_2$.
By Lemmas~\ref{res:route-3-21} and~\ref{res:route-3-29}, the two-vertex members of $\mathcal A$
are precisely the pairs $\{u,z\}$ for which $u$ is an owner of $z$.
Indeed, Lemma~\ref{res:route-3-29} ensures that this is the only critical set with its root and designated neighbour,
so it is the corresponding selected set $S_{u,v}$.
Therefore $N=4E+N_2+M$.

Put $F=|\mathcal F|$, $\xi=\sum_{c\in C}(d_{T(G)}(c)-1)$ and $\tau=|I|$.
Let $r_0$ be the number of edges in $I$ neither of whose end vertices is incident with an edge of $\mathcal F$.
Write $b_D,b_P,b_C$ for the numbers of large sets $S_{u,v}$ with $v$ in $D,P,C$, respectively.
Then $\ell=b_D+b_P+b_C$, where $\ell$ is the number of large sets in the family selected in Subsection~3.3.
For each of the $b_P$ sets, Lemma~\ref{res:route-3-37}
gives two large nontrivial critical sets with roots $c_i$ and designated neighbours $w_i$ for $i=1,2$,
where the roots $c_i$ are distinct vertices of $C_{PP}$ having owners.
For a large set $S_{u,v}$ with $v\in C$, Lemma~\ref{res:route-3-28}
gives a vertex $z$ owned by $u$ such that $vz$ is thin.
If $u$ were a centre of a subgraph given by Lemma~\ref{res:route-3-34}, the proof of that lemma would give
a two-vertex maximal nontrivial critical set with root $u$ and designated neighbour $v$.
By Lemma~\ref{res:route-3-29}, it would be the only such set, contrary to $S_{u,v}$ being large.
Thus $u$ is not a centre. If $d_{T(G)}(z)=1$, then $z\in C_{PP}$ by Lemma~\ref{lem:cubic-owners}(i),
and uniqueness of the owner would make $u$ a centre, a contradiction.
Hence $d_{T(G)}(z)=2$.
These vertices $z$ are distinct, since owners are unique and each root has at most one large set $S_{u,v}$.
Thus $b_C\le N_2$.
Let $w$ be the neighbour of $z$ outside $N(u)$.
By Lemma~\ref{res:route-3-35}, $S_{u,v}\cup\{z\}$ is a large nontrivial critical set
with root $z$ and designated neighbour $w$.

The vertices $z$ lie outside $C_{PP}$, so they and the roots $c_i$ are all
distinct. The designated neighbours $w_i$ and $w$ lie in $D\cup C$, by
Lemma~\ref{res:route-3-17} at the corners $c_i$ and Lemma~\ref{lem:cubic-owners}(iii) at the vertices $z$. Let $\ell_D,\ell_C$ count these roots according as their
designated neighbours belong to $D$ or $C$. Then $\ell_D+\ell_C=2b_P+b_C$.
The $\ell_C$ roots give distinct vertex-edge pairs $(c,cv)$, where
$cv$ is a nonthin edge of $G[C]$ and $c$ has an owner and lies outside $V(I)$.

Sum the inequality in Lemma~\ref{lem:cubic-small-set} over the $M$ sets $S\in\mathcal A$ with $|S|\ge3$.
Each vertex of $C$ belongs to at most one such set, so the sum of the
$\omega(S\cap C)$ is at most $2\sum_{e\in I}w(e)
=2\sum_{c\in V(I)}d_{\mathcal F}(c)-2\tau+2r_0$.
Here $2r_0$ accounts for the edges neither of whose end vertices is incident with an edge of $\mathcal F$.
The third term in Lemma~\ref{lem:cubic-small-set} counts twice the vertex-edge pairs $(c,f)$ with
$c\in(S\cap C)\setminus V(I)$ and $f\in\mathcal F$ incident with $c$.
No such vertex $c$ has an owner. Neither these pairs nor those at $V(I)$
include the $\ell_C$ pairs defined above. Since $\mathcal F$ has $2F$
vertex-edge pairs, summing the $\omega(S\cap C)$ term and the third term over these $M$ sets gives at most
$2(2F-\ell_C)-2\tau+2r_0$.

As shown in the proof of Theorem~\ref{res:route-4-2}, each vertex of $C$ contributes to
$\alpha(S\cap C)$ for at most two members $S$ of $\mathcal A$ with $|S|\ge3$,
once through membership in $S$ and once through membership in $N_{T(G)}(S\cap C)$.
A vertex having an owner belongs to none of these sets $S$.
Thus each of the $N_2$ vertices contributes at most once, and
$\sum_{S\in\mathcal A,\ |S|\ge3}\alpha(S\cap C)\le2\xi-N_2$.
Combining the three estimates and dividing by two gives
\begin{equation}\label{eq:route-6-1}
 2M+\ell_C\le2F-\tau+r_0+2\xi-N_2.
\end{equation}

For each edge $ab$ counted by $r_0$, there is a large nontrivial critical set
with root $a$ and designated neighbour $z\in N(a)\setminus\{b\}$, or with root $b$ and
designated neighbour $z\in N(b)\setminus\{a\}$. Indeed, the two neighbours of $a$ other than $b$ and the two neighbours of $b$ other than $a$ lie in $D$.
They do not belong to $P$ by Lemma~\ref{lem:cubic-owners}(iii).
The four edges are nonthin, so the definition of $r_0$ also excludes neighbours in $C$.
Suppose all maximal nontrivial critical sets with these roots and designated neighbours
are small. Among these sets, select $S$ of minimum cardinality, with root $a$, say.
If $N(b)\cap S=\{a\}$, then $S$ would imply that $ab$ is nonthin,
by Proposition~\ref{res:route-3-1} and Lemma~\ref{res:route-3-2}.
Thus $S$ contains a neighbour $x\in D$ of $b$ other than $a$. For every nontrivial critical set $F_0$ with root $x$ and designated neighbour $b$,
the paired set has root $b$ and designated neighbour $x$ and is contained in a
small maximal set. Thus
$|F_0|\ge q-1-b_q\ge b_q+1>|S|$, so no such $F_0$ lies in $S$.
Lemma~\ref{res:route-3-26} gives a small maximal nontrivial critical set of size at most
$|S|-2$, with root $b$ and designated neighbour $x$, contradicting minimality.

Select one such large set for each of the $r_0$ edges. Its root lies in $V(I)$.
These roots are distinct, since $I$ is a matching, and differ from the roots $c_i,z$, which have owners, and the roots in $P$ of the $b_D$ large sets $S_{u,v}$.
All $b_D+\ell_D+r_0$ large sets have distinct roots outside $D$ and
designated neighbours in $D$. By Lemma~\ref{res:route-3-36}, these designated neighbours lie
in distinct components of $G[D]$. Hence $b_D+\ell_D+r_0\le k$ and
\begin{equation}\label{eq:route-6-3}
\begin{aligned}
 2\ell+2r_0&=2b_D+\ell_D+\ell_C+b_C+2r_0\\
 &\le2k+\ell_C+N_2.
\end{aligned}
\end{equation}

By Lemma~\ref{res:route-3-8}(iii), $G[D]$ is a forest with $|D|-k$ edges.
Counting the edges between $D$ and $C\cup P$ gives
$|D|+2k=3(|C|+p)-2(t+F+E+L)$.
Since $|C|=2t-\xi$ and $n=|D|+|C|+p$, we obtain
\begin{equation}\label{eq:route-6-4}
 6t-n=4\xi+2L+2F+2E+2k-4p.
\end{equation}
Let $W=\{z\in C:d_{T(G)}(z)\ge2\}$.
Every thin edge is one of the $2E$ thin edges in the subgraphs given by Lemma~\ref{res:route-3-34},
belongs to $I$, or has an end vertex in $W$.
Since $d_{T(G)}(z)\le2(d_{T(G)}(z)-1)$ for $z\in W$, it follows that
\[
 t\le2E+\tau+\sum_{z\in W}d_{T(G)}(z)\le2E+\tau+2\xi.
\]

If $p>0$, Theorem~\ref{res:route-3-19} gives $N\ge2p-\ell+1$.
Since $N=4E+N_2+M$ and $L\ge4E+2N_2$, combining
\eqref{eq:route-6-1} and \eqref{eq:route-6-3} yields
\begin{align*}
 4p&\le8E+2N_2+2M+2\ell-2\\
   &\le2L+2F+2k+2\xi-\tau-r_0-2N_2-2.
\end{align*}
Thus $6t-n\ge2\xi+2E+\tau+r_0+2N_2+2$.
Using $t\le2E+\tau+2\xi$, we obtain $5t-n\ge2N_2+r_0+2$.

If $p=0$, then $E=L=N_2=0$. The preceding choice of a large set for
each of the $r_0$ edges, together with Lemma~\ref{res:route-3-36}, gives $r_0\le k$.
Each of the other $\tau-r_0$ edges of $I$ meets
an edge of $\mathcal F$, so $\tau-r_0\le2F$.
Consequently $t\le\tau+2\xi\le2F+r_0+2\xi$, and \eqref{eq:route-6-4} gives
$5t-n\ge2\xi+2k-r_0\ge2\xi+k$. This is positive unless $\xi=k=0$.
In that case, $D=\varnothing$ and $P=\varnothing$, so every vertex belongs to $C$;
the equality $\xi=0$ implies that every vertex is incident with exactly one thin edge. Then $n=2t$ and $5t-n=3t>0$.
\end{proof}

\subsection{The Proof of Theorem~\ref{cub:thm:main}}
We begin with three lemmas needed for the proof of Theorem~\ref{cub:thm:main}.

\begin{lemma}\label{res:route-6-3}
Let $G$ be a connected cubic bipartite graph of order at least six.
Then $G$ is matching covered. Moreover, $G$ is a brace if and only if every nontrivial cut has at least four edges.
\end{lemma}
\begin{proof}
Let $(A,B)$ be the bipartition of $G$. Hall's theorem gives a perfect matching $M$.
The edges outside $M$ form disjoint even cycles.
Selecting alternate edges of each cycle gives two further perfect matchings, which together with $M$ partition $E(G)$.
Thus $G$ is matching covered.

Suppose that $G$ is a brace. By Theorem~\ref{thm:brace-connectivity}, every cut has at least three edges.
Let $\partial(X)$ have three edges, and let $r,s$ be the numbers of these edges incident with
$A\cap X,B\cap X$, respectively.
Counting degrees in these two sets gives $r+s=3$ and $r-s=3(|A\cap X|-|B\cap X|)$.
Thus $\{r,s\}=\{0,3\}$.
The larger of $A\cap X,B\cap X$ has one more vertex than the smaller, and every edge of $\partial(X)$
is incident with a vertex of the larger set.
By Proposition~\ref{prop:bipartite-tight}, $\partial(X)$ is tight and hence trivial.
Conversely, suppose that every nontrivial cut has at least four edges.
If $\partial(X)$ is tight, its edges are incident with the majority part of $X$,
which has one more vertex than the minority part by Proposition~\ref{prop:bipartite-tight}.
Counting degrees in these parts gives $|\partial(X)|=3$.
Therefore every tight cut is trivial, and $G$ is a brace.
\end{proof}

\begin{lemma}\label{res:route-6-4}
Let $G$ be a cubic brace, and let $\partial(X)$ be a cut of four pairwise disjoint edges.
If both shores have at least six vertices, every edge of $\partial(X)$ is nonthin.
\end{lemma}
\begin{proof}
Let $(A,B)$ be the bipartition of $G$, put $q=|A|=|B|$, and let $Y=V(G)\setminus X$.
Write $r,s$ for the numbers of cut edges incident with $A\cap X,B\cap X$, respectively.
Then $r+s=4$ and $r-s=3(|A\cap X|-|B\cap X|)$, so $r=s=2$ and $|A\cap X|=|B\cap X|$.
It follows that $|A\cap Y|=|B\cap Y|$ as well.
Let $ab$ be any edge of the cut. Interchanging $A,B$ if necessary, assume $a\in A\cap X$ and $b\in B\cap Y$.
Let $a'$ be the end vertex in $A\cap X$ of the other cut edge joining $A\cap X$ to $B\cap Y$.
Since the cut is a matching, $a'\ne a$.
The sets $(A\cap X)\setminus\{a'\}$ and $B\cap Y$ each have at least two vertices,
their sizes sum to $q-1$, and $ab$ is the only edge joining them.
By Proposition~\ref{res:route-3-1}, $ab$ is nonthin.
Since $ab$ was arbitrary, every edge of the cut is nonthin.
\end{proof}

\begin{lemma}\label{res:route-6-5}
Let $G'$ be a bipartite graph obtained from two vertex-disjoint simple cubic braces
by deleting a matching of size two from each and joining the four end vertices of one deleted matching
to the four end vertices of the other by a matching of four edges. Then the following statements hold.
\begin{enumerate}[label=\textup{(\roman*)},leftmargin=*]
\item The graph $G'$ is a simple cubic brace.
\item All four added edges are nonthin in $G'$.
\end{enumerate}
\end{lemma}
\begin{proof}
Denote the two braces by $G_1,G_2$ and the deleted matchings by $M_1,M_2$.
For $i=1,2$, put $H_i=G_i-M_i$.
By Lemma~\ref{res:route-6-3}, every cut of $G_i$ has at least three edges, so $H_i$ is connected.
Suppose that $e$ is a bridge of $H_i$, and let $Y$ be the vertex set of a component of $H_i-e$.
Then $\partial_{G_i}(Y)\subseteq M_i\cup\{e\}$, so this cut has at most three edges.
By Lemma~\ref{res:route-6-3}, the cut is trivial and hence has exactly three edges, since $G_i$ is cubic.
Both edges of $M_i$ would then be incident with the vertex in the singleton shore,
contrary to the fact that $M_i$ is a matching. Thus $H_1$ and $H_2$ have no bridges.

The graph $G'$ is simple, cubic, bipartite and connected.
Let $\partial_{G'}(X)$ be a cut with at most three edges.
If both $X$ and its complement met each $V(H_i)$, each $H_i$ would have at least two edges in $\partial_{G'}(X)$,
since it is connected and has no bridge. The cut would then have at least four edges, a contradiction.
Also, the cut cannot have shores $V(H_1),V(H_2)$, since four edges join them.
After interchanging the shores and the subscripts, we may therefore assume
$\varnothing\ne X\subsetneq V(H_1)$.
Let $m$ be the number of edges of $M_1$ with both end vertices in $X$.
Each edge of $M_1\cap\partial_{G_1}(X)$ is deleted, and its end vertex in $X$ is joined to $V(G_2)$ by one new edge.
For each edge of $M_1$ with both end vertices in $X$, two new edges join its end vertices to $V(G_2)$.
Consequently $|\partial_{G'}(X)|=|\partial_{G_1}(X)|+2m$.
By Lemma~\ref{res:route-6-3}, $\partial_{G_1}(X)$ is trivial.
If $V(G_1)\setminus X=\{z\}$, at most one edge of $M_1$ could be incident with $z$, since $M_1$ is a matching.
Hence at least one edge of $M_1$ would have both end vertices in $X$.
This would give $m\ge1$ and $|\partial_{G'}(X)|\ge5$.
Thus $X$ is a singleton. Every nontrivial cut of $G'$ therefore has at least four edges,
so $G'$ is a brace by Lemma~\ref{res:route-6-3}. This proves (i).

The four added edges form the cut $\partial_{G'}(V(G_1))$ and are pairwise nonadjacent.
Its shores $V(G_1),V(G_2)$ each have at least six vertices, so (ii) follows from Lemma~\ref{res:route-6-4}.
\end{proof}

We now prove Theorem~\ref{cub:thm:main}.

\begin{proof}[\normalfont\bfseries Proof of Theorem~\ref{cub:thm:main}]
Let $G$ be a cubic brace of order $n\ge6$, and put $t=|T(G)|$.
By Theorem~\ref{thm:brace-connectivity}, $G$ is simple.
If $n\ge12$, Lemma~\ref{res:route-6-2} gives $5t>n$.
Since $t$ is an integer, $t\ge\lceil\frac{n+1}{5}\rceil$.
For $n=6,8$, the same bound follows because every brace of order at least six has at least two thin edges~\cite{CLM}.

Suppose that $n=10$. The complement of $G$ in $K_{5,5}$ is a simple $2$-regular bipartite graph,
so it is either $C_{10}$ or $C_4\cup C_6$.
If the complement is $C_4\cup C_6$, let $X$ consist of the two vertices in one part of the $4$-cycle.
Then $|N_G(X)|=3$, contrary to Theorem~\ref{thm:brace-hall}.
Thus the complement is $C_{10}$.
The rotation that moves each vertex two positions along $C_{10}$ in the same direction
is an automorphism of $G$ preserving its bipartition. It also preserves thinness.
Select a thin edge, which exists by~\cite{CLM}.
This edge and those obtained by four successive rotations are five distinct thin edges,
since their end vertices in either part are distinct. Hence $t\ge5$.

It remains to show that the coefficient $\frac15$ is best possible.
Let $k$ be a nonnegative integer.
Let $G_0$ be a tree with $k$ vertices of degree three, $k+2$ leaves, and no other vertices.
Such a tree can be obtained from a single edge by repeatedly adding two new leaves adjacent to a leaf.
For each vertex $r$ of $G_0$, we specify a graph $H_r$ and one matching for each edge of $G_0$ incident with $r$.
Take these graphs $H_r$ to be vertex-disjoint.

If $r$ is cubic in $G_0$, take $H_r$ to be the cube with bipartition
$(\{u,a_1,a_2,a_3\},\{v,b_1,b_2,b_3\})$ and edges $ub_i,va_i$ for $i=1,2,3$,
together with $a_ib_j$ for $i\ne j$.
This is $K_{4,4}$ with a perfect matching deleted.
It is matching covered by Lemma~\ref{res:route-6-3}.
Each vertex has three neighbours and each pair of vertices in the same part has four neighbours,
so it is a brace by Theorem~\ref{thm:brace-hall}.
Index the three incident edges of $G_0$ by $1,2,3$.
For each $i$, assign the matching $\{a_jb_\ell,a_\ell b_j\}$ to the edge with index $i$,
where $\{i,j,\ell\}=\{1,2,3\}$.
These three matchings partition the six edges $a_ib_j$ with $i\ne j$.

If $r$ is a leaf of $G_0$, take $H_r$ to be the hexagonal prism consisting of the cycles
$(a_1,y_1,c_1,d_2,x_2,b_2)$ and $(b_1,x_1,d_1,c_2,y_2,a_2)$
and the matching consisting of the edges $a_ib_i,x_iy_i,c_id_i$ for $i=1,2$.
One part consists of the vertices $a_i,x_i,c_i$, and the other of the vertices $b_i,y_i,d_i$, for $i=1,2$.
A cut with at most three edges cannot have both shores meeting the vertex sets of both cycles,
since at least two edges of each cycle would cross the cut.
Nor can its shores be the vertex sets of the two cycles, since six edges join them.
Thus it has a shore $X$ properly contained in the vertex set of one of the cycles.
At least two cycle edges and all $|X|$ matching edges incident with $X$ cross the cut.
Hence $2+|X|\le3$, so $X$ is a singleton.
By Lemma~\ref{res:route-6-3}, the prism is a brace.
Assign the matching $\{b_1a_2,d_1c_2\}$ to the edge of $G_0$ incident with this leaf.

For each edge $rs$ of $G_0$, delete the matching assigned to $rs$ in $H_r$ and the matching assigned to $rs$ in $H_s$,
and join their end vertices as in Lemma~\ref{res:route-6-5},
with each new edge joining opposite parts of the bipartition.
At each step the two matchings lie in distinct components, because $G_0$ is a tree.
The matchings assigned to edges of $G_0$ not yet considered remain present, since the three matchings in each cube are edge-disjoint.
By Lemma~\ref{res:route-6-5}, the component obtained at each step is a simple cubic brace.
Denote the final graph by $G'$, with bipartition $(A,B)$, and put $q=|A|=|B|$.
Since $G_0$ has at least two leaves and each corresponding graph $H_r$ has twelve vertices,
we have $|V(G')|\ge24$ and hence $q\ge12$.

For a vertex $r$ of degree three in $G_0$, the graph $G'[V(H_r)]$ consists of
the edges $ub_i,va_i$ for $i=1,2,3$.
For a leaf $r$, the graph $G'[V(H_r)]$ has edges
$a_ib_i,a_iy_i,x_ib_i,x_iy_i,x_id_i,c_iy_i,c_id_i$ for $i=1,2$, together with $a_1b_2,c_1d_2$.
Each set $V(H_r)$ has eight or twelve vertices, equally divided between $A$ and $B$.
Deleting an edge of $G_0$ gives two components.
The unions of the sets $V(H_r)$ over these components form the shores of a cut of $G'$ consisting of four pairwise nonadjacent edges.
Both shores have at least eight vertices, so all four edges of this cut are nonthin by Lemma~\ref{res:route-6-4}.
Thus every edge added between two of the graphs $H_r$ is nonthin in $G'$.

Consider a vertex $r$ of degree three in $G_0$ and an incident edge with index $i$.
Delete this tree edge, and let $V_i$ be the union of the sets $V(H_s)$ for vertices $s$
in the component not containing $r$.
Interchanging $A$ and $B$ if necessary, assume that $u\in A$.
Write $h=|A\cap V_i|=|B\cap V_i|$, and put $S=\{u\}\cup(A\cap V_i)$.
Since the four edges leaving $V_i$ are pairwise disjoint, every vertex of $G'[V_i]$ has degree at least two.
Thus every vertex of $B\cap V_i$ has a neighbour in $A\cap V_i$.
The only neighbours of $A\cap V_i$ outside $V_i$ are $b_j,b_\ell$,
where $\{i,j,\ell\}=\{1,2,3\}$.
It follows that $N_{G'}(S)=(B\cap V_i)\cup\{b_1,b_2,b_3\}$ and $|N_{G'}(S)|=h+3=|S|+2$.
Moreover, $N_{G'}(b_i)\cap S=\{u\}$.
The set $S$ has at least two vertices, and $a_1,a_2,a_3\notin S$ gives $|S|\le q-3$.
Thus $S$ is a nontrivial critical set with root $u$ and designated neighbour $b_i$.
By Proposition~\ref{res:route-3-1} and Lemma~\ref{res:route-3-2}, $ub_i$ is nonthin.
Interchanging the parts of the bipartition gives the same conclusion for $va_i$.

Now consider a leaf $r$ of $G_0$.
For $i=1,2$, the sets $\{a_i,x_i\}$ and $\{x_i,c_i\}$ are nontrivial critical sets with root $x_i$ and respective designated neighbours $d_i,b_i$;
$\{b_i,y_i\}$ and $\{y_i,d_i\}$ have root $y_i$ and respective designated neighbours $c_i,a_i$.
Indeed, the vertices of each pair have two common neighbours and distinct third neighbours, giving four neighbours in all.
The pairs $\{a_1,x_1\}$ and $\{x_1,c_1\}$ also have respective roots $a_1,c_1$
and respective designated neighbours $b_2,d_2$.
By Proposition~\ref{res:route-3-1} and Lemma~\ref{res:route-3-2}, the edges
$x_id_i,x_ib_i,y_ic_i,y_ia_i$ for $i=1,2$, and $a_1b_2,c_1d_2$ are nonthin.

Finally, put $S_1=\{y_1,b_2,y_2,d_2\}$ and $S_2=\{x_2,a_1,x_1,c_1\}$.
Their neighbourhoods in $G'$ are $\{a_1,x_1,c_1,a_2,x_2,c_2\}$ and $\{b_1,y_1,d_1,b_2,y_2,d_2\}$, respectively.
Also, $N_{G'}(x_1)\cap S_1=\{y_1\}$ and $N_{G'}(y_2)\cap S_2=\{x_2\}$.
Thus $S_1,S_2$ are nontrivial critical sets with respective roots $y_1,x_2$
and respective designated neighbours $x_1,y_2$.
By Proposition~\ref{res:route-3-1} and Lemma~\ref{res:route-3-2}, $x_1y_1,x_2y_2$ are nonthin.
Only the four edges $a_ib_i,c_id_i$, for $i=1,2$, of $G'[V(H_r)]$ can therefore be thin in $G'$.
No vertices are deleted in obtaining $G'$, so $|V(G')|=8k+12(k+2)=20k+24$.
Every thin edge of $G'$ lies in $G'[V(H_r)]$ for some leaf $r$,
and each of these $k+2$ subgraphs contains at most four thin edges. Hence $|T(G')|\le4(k+2)$.
The lower bound already proved therefore gives
\[
 \frac15<\frac{|T(G')|}{|V(G')|}\le\frac{4k+8}{20k+24}.
\]
The right-hand side tends to $\frac15$ as $k$ tends to infinity.
Thus the infimum in Theorem~\ref{cub:thm:main} is $\frac15$.
\end{proof}




\end{document}